# BATCH MEANS AND SPECTRAL VARIANCE ESTIMATORS IN MARKOV CHAIN MONTE CARLO

By James M. Flegal and Galin L. Jones[1]

*University of California, Riverside and University of Minnesota*

Calculating a Monte Carlo standard error (MCSE) is an important step in the statistical analysis of the simulation output obtained from a Markov chain Monte Carlo experiment. An MCSE is usually based on an estimate of the variance of the asymptotic normal distribution. We consider spectral and batch means methods for estimating this variance. In particular, we establish conditions which guarantee that these estimators are strongly consistent as the simulation effort increases. In addition, for the batch means and overlapping batch means methods we establish conditions ensuring consistency in the mean-square sense which in turn allows us to calculate the optimal batch size up to a constant of proportionality. Finally, we examine the empirical finite-sample properties of spectral variance and batch means estimators and provide recommendations for practitioners.

**1. Introduction.** Suppose $\pi$ is a probability distribution with support $\mathsf{X}$ and the goal is to calculate $E_\pi g := \int_{\mathsf{X}} g(x)\pi(dx)$ where $g$ is a real-valued, $\pi$-integrable function. In many situations, $\pi$ is sufficiently complex so that we may have to rely on Markov chain Monte Carlo (MCMC) methods to estimate $E_\pi g$. As is now widely recognized (Liu [30], Robert and Casella [41]), we can often simulate a Harris ergodic (i.e., aperiodic, $\pi$-irreducible, positive Harris recurrent) Markov chain on $\mathsf{X}$ having invariant distribution $\pi$ and easily estimate $E_\pi g$. Specifically, letting $X = \{X_1, X_2, X_3, \dots\}$ denote such a Markov chain, with probability one and for any initial distribution

$$\bar{g}_n := \frac{1}{n}\sum_{i=1}^{n} g(X_i) \to E_\pi g \qquad \text{as } n \to \infty. \tag{1}$$

Received November 2008; revised April 2009.
[1]Supported in part by NSF Grant DMS-08-06178.

The approximate sampling distribution of the Monte Carlo error, $\bar{g}_n - E_\pi g$, is available via a Markov chain central limit theorem (CLT) if there exists a constant $\sigma_g^2 \in (0, \infty)$ such that

$$\sqrt{n}(\bar{g}_n - E_\pi g) \xrightarrow{d} \mathrm{N}(0, \sigma_g^2) \qquad \text{as } n \to \infty. \tag{2}$$

The conditions we will require in our theoretical work below are sufficient to guarantee (2) for any initial distribution. In fact, our conditions will imply the stronger Markov chain functional central limit theorem. See Jones [23] and Roberts and Rosenthal [44] for broader discussion of the conditions for a Markov chain CLT.

Obtaining a good estimate of $\sigma_g^2$, say $\hat{\sigma}_n^2$, is an important step in the statistical analysis of the observed sample path for at least two reasons: (1) it can be used to construct asymptotically valid confidence intervals for $E_\pi g$ and (2) it is a key component of rigorous rules for deciding when to terminate the simulation. We describe these approaches more fully below but the interested reader is directed to Flegal, Haran and Jones [12] and Jones et al. [24] for more detail and comparisons with other methods. In particular, these papers demonstrate that terminating the simulation based on $\hat{\sigma}_n^2$ is superior to the common practice of terminating based on convergence diagnostics.

The simplest approach to stopping an MCMC experiment is a *fixed-time* rule. Specifically, the Markov chain simulation is run for a predetermined number of iterations, using $\bar{g}_n$ to estimate $E_\pi g$. If $\hat{\sigma}_n^2$ is a consistent estimator of $\sigma_g^2$, a valid Monte Carlo standard error (MCSE) of $\bar{g}_n$ is given by $\hat{\sigma}_n/\sqrt{n}$ and an asymptotically valid interval estimator of $E_\pi g$ is given in the usual way

$$\bar{g}_n \pm t_* \frac{\hat{\sigma}_n}{\sqrt{n}},$$

where $t_*$ is an appropriate Student's $t$ quantile. Reporting an interval estimate of $E_\pi g$, or at least the MCSE, will allow independent evaluation of the quality of the reported point estimate. Of course, it is possible the interval estimate is undesirably wide and hence the simulation should continue. This naturally leads to a sequential approach to terminating the simulation the first time the interval is sufficiently narrow, i.e., a *fixed-width* rule, where the total simulation effort is random. Formally, the user specifies a desired half-width $\epsilon$ and the simulation terminates the first time the following inequality is satisfied

$$t_* \frac{\hat{\sigma}_n}{\sqrt{n}} + p(n) < \epsilon, \tag{3}$$

where $p(n)$ is a positive function such that as $n \to \infty$, $p(n) = o(n^{-1/2})$. Glynn and Whitt [17] established that if a functional central limit theorem holds and $\hat{\sigma}_n^2 \to \sigma_g^2$ with probability 1 as $n \to \infty$, then the interval at (3) is asymptotically valid in the sense that the desired coverage probability is obtained as $\epsilon \to 0$. Moreover, the role of $p(n)$ is to ensure that for small values of $\epsilon$ the simulation effort $n$ is large. Letting $n^*$ be the desired minimum simulation effort, an often useful choice is $p(n) = \epsilon I(n \leq n^*) + n^{-1}$ where $I(\cdot)$ is the usual indicator function.

Due to the inherent serial correlation in the Markov chain, generally $\sigma_g^2 \neq \mathrm{Var}_\pi\, g$, and hence estimating $\sigma_g^2$ requires specialized techniques such as non-overlapping batch means (BM), overlapping batch means (OBM), spectral variance (SV) methods and regenerative simulation (RS). We study asymptotic, specifically strong consistency and mean-square consistency, and finite-sample properties of the BM, OBM and SV procedures. Strong consistency of the BM and RS estimators was addressed by Jones et al. [24] and Hobert et al. [20], respectively. In the current work we develop conditions for the strong consistency of OBM and SV procedures. More specifically, we require that the Markov chain be geometrically ergodic while existing results on the consistency of OBM and SV methods require the chain to be uniformly ergodic; for a definition of geometric and uniform ergodicity see Tierney [50]. Consistency of RS and BM also require geometric ergodicity, however, SV methods require a slightly stronger moment condition on $g$ compared to RS, BM and OBM. Overall, as we discuss in Sections 2 and 3, our results significantly weaken the existing regularity conditions guaranteeing strong consistency for SV and OBM. It is also worth emphasizing that the results on strong consistency of BM, OBM, RS and SV do not require a stationary Markov chain and hence from a theoretical perspective burn-in is not required. Of course, the initial distribution should be carefully chosen since it will impact the finite sample properties of the chain and the resulting output analysis.

Establishing that a given Markov chain is geometrically or uniformly ergodic can be challenging. On the other hand, there has been a substantial amount of effort directed towards doing just this in the context of MCMC. For example, we know that Metropolis–Hastings samplers with state-independent proposals can be uniformly ergodic (Tierney [50]) but such situations are uncommon in realistic settings. Standard random walk Metropolis–Hastings chains on $\mathbb{R}^d$, $d \geq 1$ cannot be uniformly ergodic but may still be geometrically ergodic (see Mengersen and Tweedie [35]). An incomplete list of other research on establishing convergence rates of Markov chains used in MCMC is given by Geyer [13], Jarner and Hansen [21], Meyn and Tweedie [37] and Neath and Jones [39] who considered Metropolis–Hastings algorithms and Hobert and Geyer [19], Hobert et al. [20], Johnson and Jones [22], Jones and Hobert [26], Marchev and Hobert [33], Roberts and Polson [42], Roberts and Rosenthal [43], Rosenthal [45, 46], Roy and Hobert [47], Tan and Hobert [49] and Tierney [50] who examined Gibbs samplers.

Optimal batch size selection for BM and OBM is a long-standing open problem. Song and Schmeiser [48] propose an approach that minimizes the asymptotic mean-squared error. Thus we also prove the BM and OBM estimators are *mean-square consistent*, or

$$\text{MSE}(\hat{\sigma}_n^2) := E_\pi(\hat{\sigma}_n^2 - \sigma_g^2)^2 \to 0 \qquad \text{as } n \to \infty. \tag{4}$$

Our work on establishing (4) allows us to argue that the asymptotically optimal batch size in terms of MSE is proportional to $n^{1/3}$. This is similar to the conclusions of others (Chien, Goldsman and Melamed [6], Song and Schmeiser [48]); however, our mixing and moment conditions are much weaker. Specifically, the previous work requires a uniformly ergodic Markov chain and an absolute twelfth moment where we require only geometric ergodicity and a bit more than a fourth moment.

The BM, OBM, RS and SV procedures are all easy to implement. RS is sometimes viewed as the standard by which others should be judged (see, e.g., Bratley, Fox and Schrage [4], Jones and Hobert [25]). However, RS may require additional theoretical work which could be an obstacle to some, and moreover, has been found to be problematic in high-dimensional settings (Gilks, Roberts and Sahu [14]) and in variable-at-a-time Metropolis–Hastings settings (Neath and Jones [39]). Of the alternatives to RS, OBM and SV procedures have the reputation of being more efficient than BM. For example, results in Section 3.3 show the ratio of the variance of the OBM estimator to the variance of the BM estimator converges to 2/3 as the simulation effort increases. However, asymptotics alone do not appear to give a clear picture of which method should be preferred. Hence we investigate the empirical finite-sample properties of these methods in the context of two examples. Discussion of our findings is given in Section 4.3 but, on balance, some of the SV methods appear to be superior.

The rest of this paper is organized as follows. In Section 2 we establish strong consistency of SV procedures. Next, Section 3 contains a discussion of the asymptotic properties of BM and OBM procedures. The finite sample properties of the various methods are studied in two examples in Section 4 where we also give some general recommendations. Finally, most of the technical details and proofs are presented in appendices.

**2. Spectral estimation.** In this section, we define a class of estimators of $\sigma_g^2$ and establish conditions which guarantee their strong consistency. First, define the lag $s$ autocovariance $\gamma(s) = \gamma(-s) := E_\pi[Y_t Y_{t+s}]$ where $Y_i := g(X_i) - E_\pi g$ for $i = 1, 2, 3, \ldots$. Consider estimating $\gamma(s)$ with

$$\gamma_n(s) = \gamma_n(-s) := n^{-1} \sum_{t=1}^{n-s} (Y_t - \bar{Y}_n)(Y_{t+s} - \bar{Y}_n),$$

where $\bar{Y}_n := n^{-1}\sum_{i=1}^n Y_i$. If $E_\pi g^2 < \infty$, then for fixed $s$, $\gamma_n(s) \to \gamma(s)$ w.p. 1 as $n \to \infty$.

One could use the sum of the $\gamma_n(s)$ to estimate $\sigma_g^2$, though this turns out to be a poor estimator (see Anderson [1] and Bratley, Fox and Schrage [4]). Instead we will investigate a truncated and weighted estimator version called the *spectral variance estimator*

$$\hat{\sigma}_S^2 := \sum_{s=-(b_n-1)}^{b_n-1} w_n(s)\gamma_n(s),$$

where $w_n(\cdot)$ is the *lag window*, and $b_n$ is the *truncation point*.

Our theoretical results require the following assumptions on the lag window, truncation point and one-step Markov transition kernel associated with $X$, denoted $P(x, A)$.

Assumption 1. The lag window $w_n(\cdot)$ is an even function defined on the integers such that

$$\begin{aligned} |w_n(s)| &\le 1 && \text{for all } n \text{ and } s, \\ w_n(0) &= 1 && \text{for all } n, \\ w_n(s) &= 0 && \text{for all } |s| \ge b_n. \end{aligned}$$

Assumption 2. Let $b_n$ be an integer sequence such that $b_n \to \infty$ and $n/b_n \to \infty$ as $n \to \infty$ where $b_n$ and $n/b_n$ are monotonically nondecreasing.

Assumption 3. There exists a function $s : \mathsf{X} \to [0,1]$ and a probability measure $Q$ such that $P(x,\cdot) \ge s(x)Q(\cdot)$ for all $x \in \mathsf{X}$.

The main result of this section follows.

Theorem 1. *Let $X$ be a geometrically ergodic Markov chain with invariant distribution $\pi$ and $g : \mathsf{X} \to \mathbb{R}$ be a Borel function with $E_\pi|g|^{4+\delta+\epsilon} < \infty$ for some $\delta > 0$ and $\epsilon > 0$. Suppose Assumptions 1, 2 and 3 hold and define $\Delta_1 w_n(k) = w_n(k-1) - w_n(k)$ and $\Delta_2 w_n(k) = w_n(k-1) - 2w_n(k) + w_n(k+1)$. Further suppose* (a) *$b_n n^{-1}\sum_{k=1}^{b_n} k|\Delta_1 w_n(k)| \to 0$ as $n \to \infty$;* (b) *there exists a constant $c \ge 1$ such that $\sum_n (b_n/n)^c < \infty$;* (c) *$b_n n^{-1}\log n \to 0$ as $n \to \infty$;* (d)

$$b_n n^{2\alpha}(\log n)^3\left(\sum_{k=1}^{b_n} |\Delta_2 w_n(k)|\right)^2 \to 0 \qquad \text{as } n \to \infty,$$

*and*

$$n^{2\alpha}(\log n)^2 \sum_{k=1}^{b_n} |\Delta_2 w_n(k)| \to 0 \qquad \text{as } n \to \infty,$$

*where* $\alpha = 1/(4+\delta)$*; and* (e) $b_n^{-1} n^{2\alpha} \log n \to 0$ *as* $n \to \infty$*. Then for any initial distribution, with probability 1,* $\hat{\sigma}_S^2 \to \sigma_g^2$ *as* $n \to \infty$*.*

PROOF. See Appendix B.1. □

REMARK 1. It is convenient in applications to take $b_n = \lfloor n^\nu \rfloor$ for some $0 < \nu < 1$ in which case conditions (b) and (c) of Theorem 1 are automatically satisfied.

REMARK 2. Assumption 3 is not critical. First, we do not require the actual value of $s$ or $Q$ at any point in this paper. Thus, unlike RS which is entirely based on assumption 3, there is no practical point in searching for a good minorization. Secondly, recall that fundamental Markov chain theory (see Chapter 5 of Meyn and Tweedie [36]) ensures the existence of an integer $n_0$ for which a minorization condition holds for the $n_0$-step kernel, that is, $P^{n_0}$. If we cannot establish the one-step minorization in Assumption 3, but we can establish an $n_0$-step minorization, then we would just use the chain with kernel $P^{n_0}$. This is reasonable since the $P^{n_0}$ inherits the stability properties of $P$.

REMARK 3. Anderson [1] gives an extensive collection of lag windows satisfying Assumption 1. It is useful to consider the applicability of the conditions (a) and (d) of Theorem 1 for some of these windows.

*Simple truncation*: set $w_n(k) = I(|k| < b_n)$, then it is easy to see that condition (d) requires $4 b_n n^{2\alpha}(\log n)^3 \to 0$ which obviously cannot hold.

*Blackman–Tukey*: let $w_n(k) = [1 - 2a + 2a\cos(\pi|k|/b_n)] I(|k| < b_n)$ where $a > 0$. When $a = 1/4$ this is the *Tukey–Hanning* window. That condition (a) holds if $b_n^2 n^{-1} \to 0$ as $n \to \infty$ while condition (d) is satisfied if $b_n^{-1} n^{2\alpha}(\log n)^3 \to 0$ as $n \to \infty$ follows easily from Lemma 7 in Appendix B.2.

*Parzen*: let $w_n(k) = [1 - |k|^q / b_n^q] I(|k| < b_n)$ for $q \in \mathbb{Z}_+ = \{1, 2, 3, \ldots\}$. When $q = 1$ this window is the *modified Bartlett* window and deserves to be singled out because of its connection to the method of overlapping batch means which we will consider later. In this case $\Delta_1 w_n(k) = \Delta_2 w_n(b_n) = b_n^{-1}$ for $k = 1, \ldots, b_n$ so that condition (a) requires $b_n^2 n^{-1} \to 0$ as $n \to \infty$. Next, note that $\Delta_2 w_n(k) = 0$ for $k = 1, 2, \ldots, b_n - 1$. Thus condition (d) is satisfied if $b_n^{-1} n^{2\alpha}(\log n)^3 \to 0$ as $n \to \infty$.

When $q \geq 2$ it is easy to show that the conditions of Lemma 7 in Appendix B.2 hold and hence that conditions (a) and (d) will be satisfied under exactly the same conditions as required for the modified Bartlett lag window. That is, if $b_n^2 n^{-1} \to 0$ and $b_n^{-1} n^{2\alpha} (\log n)^3 \to 0$ as $n \to \infty$.

*Scale-parameter modified bartlett*: let $w_n(k) = [1 - \lambda|k|/b_n] I(|k| < b_n)$ where $\lambda \neq 1$ is a positive constant. Then $\Delta_1 w_n(k) = \lambda b_n^{-1}$ for $k = 1, \dots, b_n - 1$ and $\Delta_1 w_n(b_n) = 1 - \lambda + \lambda b_n^{-1}$ so that condition (a) becomes $(1-\lambda) b_n n^{-1} + \lambda b_n (b_n + 1)(2n)^{-1} \to 0$, that is, $b_n^2 n^{-1} \to 0$ as $n \to \infty$. On the other hand, there is trouble with condition (d). Note that $\Delta_2 w_n(b_n) = 1 - \lambda + \lambda b_n^{-1}$ and $\Delta_2 w_n(b_n - 1) = -1 + \lambda$ but $\Delta_2 w_n(k) = 0$ for $k = 1, 2, \dots, b_n - 2$. Hence, as $n \to \infty$, $\sum_{k=1}^{b_n} |\Delta_2 w_n(k)|$ does not converge to 0. Thus condition (d) cannot hold.

REMARK 4. Damerdji [9, 10] has previously addressed strong consistency of $\sigma_S^2$. However, our result substantially weakens the regularity conditions for Harris ergodic Markov chains. In particular, Damerdji's approach requires a uniformly ergodic chain whereas Theorem 1 requires only geometric ergodicity. Also, instead of condition (d), Damerdji's result requires

$$b_n n^{1-2\alpha'} (\log n) \left( \sum_{k=1}^{b_n} |\Delta_2 w_n(k)| \right)^2 \to 0 \qquad \text{as } n \to \infty, \quad \text{and}$$

(5)

$$n^{1-2\alpha'} \sum_{k=1}^{b_n} |\Delta_2 w_n(k)| \to 0 \qquad \text{as } n \to \infty,$$

where $0 < \alpha' \leq (\delta - 2 + \epsilon)/(24 + 12(\delta - 2 + \epsilon))$. This requirement is not particularly useful when $b_n = \lfloor n^\nu \rfloor$. For example, consider using the modified Bartlett lag window. Then just as in Theorem 1, Damerdji requires $n^{2\nu - 1} \to 0$ as $n \to \infty$ but (5) requires $n^{1-2\alpha'-\nu} (\log n) \to 0$ as $n \to \infty$ and there is no $\nu$ value that satisfies both of these requirements.

Finally, Damerdji required what we view to be an unnatural regularity condition. Specifically, for large enough $n$, $b_n^{-1} \sum_{i=n-b_n+1}^{n} Y_i^2$ is almost surely bounded above. An inspection of our proof will show that we can weaken the moment condition to $E_\pi |g|^{2+\delta+\epsilon} < \infty$ for some $\delta > 0$ and $\epsilon > 0$ by making the same assumption as Damerdji.

**3. Batch means.** In nonoverlapping batch means the output is broken into blocks of equal size. Suppose the algorithm is run for a total of $n = a_n b_n$ iterations and for $k = 0, \dots, a_n - 1$ define $\bar{Y}_k := b_n^{-1} \sum_{i=1}^{b_n} Y_{kb_n + i}$. The BM estimate of $\sigma_g^2$ is

(6)
$$\hat{\sigma}_{\mathrm{BM}}^2 = \frac{b_n}{a_n - 1} \sum_{k=0}^{a_n - 1} (\bar{Y}_k - \bar{Y}_n)^2.$$

It is well known that generally (6) is not a consistent estimator of $\sigma_g^2$ (Glynn and Iglehart [15], Glynn and Whitt [16]). On the other hand, Jones et al. [24] show that if the batch size and number of batches are allowed to increase as the overall length of the simulation does (e.g., by setting $a_n = b_n = \lfloor n^{1/2} \rfloor$), then $\hat{\sigma}^2_{\mathrm{BM}} \to \sigma_g^2$ with probability one as $n \to \infty$. However, Jones et al. [24] found that the finite sample properties can be less desirable than expected; thus we consider OBM.

OBM is a generalization of BM but it is also well known that the OBM estimator is equal, except for some end-effect terms, to the SV estimator arising from the modified Bartlett lag window—a relationship we exploit later. Note that there are $n - b_n + 1$ batches of length $b_n$ indexed by $j$ running from 0 to $n - b_n$ and define $\bar{Y}_j(b_n) = b_n^{-1} \sum_{i=1}^{b_n} Y_{j+i}$. The OBM estimator of $\sigma_g^2$ results from averaging across all batches and is defined as

$$\hat{\sigma}^2_{\mathrm{OBM}} = \frac{n b_n}{(n - b_n)(n - b_n + 1)} \sum_{j=0}^{n-b_n} (\bar{Y}_j(b_n) - \bar{Y}_n)^2. \tag{7}$$

The next result establishes strong consistency of the OBM estimator.

THEOREM 2. *Let $X$ be a geometrically ergodic Markov chain with invariant distribution $\pi$ but any initial distribution and $g : \mathsf{X} \to \mathbb{R}$ be a Borel function with $E_\pi |g|^{2+\delta+\epsilon} < \infty$ for some $\delta > 0$ and $\epsilon > 0$. Suppose Assumptions 2 and 3 hold. Further suppose* (a) *there exists a constant $c \geq 1$ such that $\sum_n (b_n/n)^c < \infty$;* (b) *$b_n n^{-1} \log n \to 0$ as $n \to \infty$;* (c) *$n^{2\alpha} (\log n)^3 / b_n \to 0$ as $n \to \infty$ where $\alpha = 1/(2+\delta)$;* (d) *there exists an integer $n_0$ and a constant $c_1$ such that for all $n \geq n_0$ we have $\log n / b_n \leq c_1$; and* (e) *$b_n^2 n^{-2} \log \log n \to 0$ and $b_n^4 n^{-3} \log \log n \to 0$ as $n \to \infty$, then as $n \to \infty$, $\hat{\sigma}^2_{\mathrm{OBM}} \to \sigma_g^2$ w.p. 1.*

For the proof, see Appendix B.3.

REMARK 5. It is possible to obtain strong consistency of $\hat{\sigma}^2_{\mathrm{OBM}}$ directly from Theorem 1 using the modified Bartlett lag window. However, Theorem 1 requires a stronger moment condition on $g$. In addition, to meet condition (a) of Theorem 1 would require that $b_n^2 n^{-1} \to 0$ as $n \to \infty$; hence the conditions of Theorem 2 are weaker.

REMARK 6. Chan and Geyer [5] established a CLT as at (2) under the same moment and mixing conditions stated in our Theorem 2. Moreover, without assuming reversibility this moment condition cannot be weakened to just a second moment for geometrically ergodic chains (see Häggström [18]).

COROLLARY 1. *Let $X$ be a geometrically ergodic Markov chain with invariant distribution $\pi$ and $g : \mathsf{X} \to \mathbb{R}$ be a Borel function with $E_\pi|g|^{2+\delta+\epsilon} < \infty$ for some $\delta > 0$ and $\epsilon > 0$. Suppose Assumption 3 holds and $b_n = \lfloor n^\nu \rfloor$ where $3/4 > \nu > (1+\delta/2)^{-1}$, then $\hat{\sigma}^2_{\mathrm{OBM}} \to \sigma^2_g$ w.p. 1.*

PROOF. This follows easily by verifying the conditions of Theorem 2. □

REMARK 7. Damerdji [10] and Jones et al. [24] show the nonoverlapping batch means estimator $\hat{\sigma}^2_{\mathrm{BM}}$ is strongly consistent for uniformly and geometrically ergodic chains, respectively, under similar conditions as those required for Theorem 2. For example, Jones et al. [24] show that if $b_n = \lfloor n^\nu \rfloor$, then strong consistency requires $1 > \nu > (1+\delta/2)^{-1} > 0$ for $\delta > 0$ resulting in weaker conditions on $\nu$ than those required in Corollary 1.

3.1. *Mean-square consistency.* We now turn our attention to showing that $\hat{\sigma}^2_n$ is consistent in the mean-square sense, that is, (4). Recall that strong consistency and mean-square consistency do not generally imply each other, hence we cannot directly appeal to the results in the previous section, although they will be useful in our proofs.

There are some existing results on this problem. For example, consider the following results from Chien, Goldsman and Melamed [6]. If $X$ is a stationary uniformly ergodic Markov chain, Assumption 2 holds, and $E_\pi g^{12} < \infty$, then as $b_n \to \infty$ (whether $a_n \to \infty$ or not)

$$b_n \operatorname{Bias}[\hat{\sigma}^2_{\mathrm{BM}}] = \Gamma + o(1), \tag{8}$$

where $\Gamma := -2\sum_{s=1}^{\infty} s\gamma(s)$ which is well defined and finite. Also, if $a_n \to \infty$ and $b_n \to \infty$

$$\frac{n}{b_n} \operatorname{Var}(\hat{\sigma}^2_{\mathrm{BM}}) = 2\sigma^4_g + o(1). \tag{9}$$

Combining (8) and (9) imply (4) for $\hat{\sigma}^2_{\mathrm{BM}}$ as $a_n \to \infty$ and $b_n \to \infty$. The next result establishes (4) where $\hat{\sigma}^2_n$ is $\hat{\sigma}^2_{\mathrm{BM}}$ or $\hat{\sigma}^2_{\mathrm{OBM}}$ under weaker mixing and moment conditions.

THEOREM 3. *Let $\hat{\sigma}^2_n$ be either the BM or OBM estimator of $\sigma^2_g$, and let $X$ be a stationary geometrically ergodic Markov chain with invariant distribution $\pi$ and $g : \mathsf{X} \to \mathbb{R}$ be a Borel function with $E_\pi|g|^{4+\delta+\epsilon} < \infty$ for some $\delta > 0$ and $\epsilon > 0$. Suppose Assumptions 2 and 3 hold and $E_\pi C^4 < \infty$ where $C$ is defined at (14). If $b_n^{-1} n^{2\alpha} (\log n)^3 \to 0$ as $n \to \infty$ where $\alpha = 1/(4+\delta)$, then $\operatorname{MSE}(\hat{\sigma}^2_n) \to 0$ as $n \to \infty$.*

REMARK 8. The proof of Theorem 3 with results in Damerdji [11] show that the conclusions also hold for uniformly ergodic chains if we replace the above condition on $b_n$ with $b_n^{-1} n^{1-2\alpha'}(\log n) \to 0$ as $n \to \infty$ where $\alpha' \leq (\delta+2+\epsilon)/(24+12(\delta+2+\epsilon))$.

REMARK 9. Results in Damerdji [11] and Philipp and Stout [40] show that when the Markov chain is uniformly ergodic

$$C(\omega) = 2 + \sum_{\tau_1(\omega)}^{\tau_2(\omega)-1} |g(X_j(\omega)|,$$

where $\tau_i$, $i = 1, 2$ are the first two times at which regenerations occur in the split-chain (see Hobert et al. [20], Jones and Hobert [25] and Mykland, Tierney and Yu [38] for an introduction to the split chain). Lemma 2 in Bednorz and Latuszyński [2] shows that, in this case, $E_\pi C^4 < \infty$. When the chain is only geometrically ergodic the representation of $C$ is not as clear but we suspect that the results of Bednorz and Latuszyński [2] can again be used to establish the moment condition on $C$ since it is still defined in terms of the split chain (see Csáki and Csörgő [7]).

3.2. *Optimal batch sizes in terms of MSE.* In this section, we will use the previous results to calculate optimal batch sizes. Chien, Goldsman and Melamed [6] and Song and Schmeiser [48] study the case of BM. Combining (8) and (9) yields

$$\mathrm{MSE}(\hat{\sigma}^2_{\mathrm{BM}}) = \frac{\Gamma^2}{b_n^2} + \frac{2 b_n \sigma_g^4}{n} + o\left(\frac{1}{b_n^2}\right) + o\left(\frac{b_n}{n}\right).$$

It is easy to use the above expression to see that $\mathrm{MSE}(\hat{\sigma}^2_{\mathrm{BM}})$ will be minimized asymptotically by selecting the optimal batch size of

$$\hat{b}^* := \left(\frac{\Gamma^2 n}{\sigma_g^4}\right)^{1/3}.$$

Notice that this optimal batch size is dependent on $\Gamma^2/\sigma_g^4$ which is typically an unknown parameter relating to the process. However, this result implies that the optimal batch size should increase proportionally to $n^{1/3}$.

The main result of this section follows.

THEOREM 4. *Let $X$ be a stationary geometrically ergodic Markov chain with invariant distribution $\pi$ and $g : \mathsf{X} \to \mathbb{R}$ be a Borel function with $E_\pi |g|^{4+\delta+\epsilon} < \infty$ for some $\delta > 0$ and $\epsilon > 0$. Suppose Assumptions 2 and 3* *hold and* $E_\pi C^4 < \infty$ *where* $C$ *is defined at* (14). *If* $b_n^{-1} n^{1/2+\alpha} (\log n)^{3/2} \to 0$ *as* $n \to \infty$ *where* $\alpha = 1/(4+\delta)$, *then*

$$\frac{n}{b_n} \operatorname{Var}(\hat{\sigma}_n^2) = c\sigma_g^4 + o(1), \tag{10}$$

*where* $c = 2$ *for BM and* $c = 4/3$ *for OBM.*

Remark 10. The proof of Theorem 4 coupled with results from Damerdji [11] can be applied to uniformly ergodic Markov chains. Specifically, the condition on $b_n$ would be changed to $b_n^{-1} n^{1-\alpha'} (\log n)^{1/2} \to 0$ as $n \to \infty$ where $\alpha' \le (\delta + 2 + \epsilon)/(24 + 12(\delta + 2 + \epsilon))$. However, if $b_n = \lfloor n^\nu \rfloor$, then $\nu \in (11/12, 1)$ whereas Theorem 4 requires $\nu \in (1/2 + \alpha, 1)$.

Combining (8) and (10) yields

$$\operatorname{MSE}(\hat{\sigma}_{\mathrm{OBM}}^2) = \frac{\Gamma^2}{b_n^2} + \frac{4 b_n \sigma_g^4}{3n} + o\left(\frac{1}{b_n^2}\right) + o\left(\frac{b_n}{n}\right),$$

which will be minimized asymptotically by selecting a batch size of

$$\hat{b}^* := \left(\frac{8\Gamma^2 n}{3\sigma_g^4}\right)^{1/3}.$$

3.3. *OBM versus BM.* Comparing the conditions of Theorem 2 with the conditions of Proposition 3 in Jones et al. [24] which addresses the strong consistency of the BM estimator, we see that conditions (d) and (e) of Theorem 2 are not required for BM. For mean-square consistency with OBM, Theorems 3 and 4 require a moment condition on $C$ which is not necessary in Chien, Goldsman and Melamed [6] and Song and Schmeiser [48]; however, the moment conditions on $g$ and mixing conditions on the Markov chain are much weaker in our results. Moreover, from an implementation point of view it is clear that OBM will require more computational resources than BM.

Why might we use OBM? Often $\hat{\sigma}_{\mathrm{OBM}}^2$ has a lower asymptotic variance compared to $\hat{\sigma}_{\mathrm{BM}}^2$. Specifically, note that (10) yields

$$\frac{\operatorname{Var}(\hat{\sigma}_{\mathrm{OBM}}^2)}{\operatorname{Var}(\hat{\sigma}_{\mathrm{BM}}^2)} \to \frac{2}{3}$$

as $n \to \infty$ (see Meketon and Schmeiser [34] for the same result under different assumptions). Also, in an effort to reduce the computational demands, Welch [52] argues that most of this benefit can be achieved by a modest amount of overlapping. For example, using a batch of size 64 and splitting the batch into 4 sub-batches, then we only need the overlapping batches (of length 64) starting at $X_1, X_{17}, X_{33}, X_{49}, X_{65}, \ldots$.

**4. Examples.** In this section, we compare BM, OBM, RS and SV via their finite sample properties in two examples; one of which is a simple AR(1) model while the other is a more realistic Bayesian probit regression model.

4.1. AR*(1) model.* Consider the following AR(1) model:

$$X_i = \rho X_{i-1} + \epsilon_i \qquad \text{for } i = 1, 2, \ldots,$$

where the $\epsilon_i$ are i.i.d. $\mathrm{N}(0,1)$. As long as $|\rho| < 1$ this Markov chain is geometrically ergodic with invariant distribution $\mathrm{N}(0, 1/(1-\rho^2))$. Also, $\mathrm{cov}_\pi(X_1, X_i) = \rho^{i-1}/(1-\rho^2)$.

Consider estimating $E_\pi X$ with $\bar{x}_n$. Clearly, a CLT holds and, indeed, $\sigma_x^2 = 1/(1-\rho)^2$. Thus, if $\hat{\sigma}_n^2$ is a strongly consistent estimate of $\sigma_x^2$, an asymptotically valid confidence interval is given by

$$\bar{x}_n \pm t_* \frac{\hat{\sigma}_n}{\sqrt{n}}, \tag{11}$$

where $t_*$ is an appropriate Student's $t$ quantile. Our goal is to evaluate the finite sample properties of this interval when $\hat{\sigma}_n$ is produced by BM, OBM and SV methods with $b_n = \lfloor n^\nu \rfloor$ for some $\nu$ and $\rho \in \{0.5, 0.95\}$. For SV estimators, we will consider the Tukey–Hanning window (TH) and the modified Bartlett window (Brt). For BM the degrees of freedom for $t_*$ are $a_n - 1$ and for OBM, TH and Brt the degrees of freedom for $t_*$ are $n - b_n$.

We compare the effects of using different batch sizes and variance estimation techniques on the coverage probabilities of a nominal 95% interval as at (11) for the two settings for $\rho$. This comparison is based upon the results of 2000 independent replications of the following procedure. In each replication we simulated the AR(1) chain for 1e5 iterations. We then constructed the interval estimate using BM, OBM, Brt and TH with three different sampling plans, $b_n = \lfloor n^\nu \rfloor$ where $\nu \in \{1/3, 1/2, 2/3\}$ at five different points (1e3, 5e3, 1e4, 5e4 and 1e5 iterations).

The results with $\rho = 0.5$ are summarized in Table 1. In the calculations with $\nu = 1/3$ and $\nu = 1/2$, all of the calculated coverage probabilities are within 2 standard errors of the nominal 0.95 level when at least 5e3 iterations are used. We can also see that for all the settings, the coverage probabilities improve as the number of iterations increase. The choice of $\nu = 2/3$ seems to slightly underestimate the coverage probabilities for small numbers of iterations. Basically, when $\rho = 0.5$ the estimation problem is relatively easy, and all the methods and settings seem to perform well.

Consider the case with $\rho = 0.95$. Table 2 shows the observed coverage probabilities. We can see the coverage probabilities get closer to the nominal 0.95 level as the number of iterations increases. Also, as $\nu$ increases, the confidence intervals become more accurate because the strong correlation in the observations is better captured with the larger batch size or truncation point. In this case, the choice of $\nu = 1/3$ performs much worse than the other options analyzed.

TABLE 1

*Table of coverage probabilities for 2000 replications using the* AR*(1) example with* $\rho = 0.5$*. All calculations were based on the nominal level of 0.95. The standard errors for these numbers are easily calculated as* $\sqrt{\hat{p}(1-\hat{p})/2000}$ *which results in a largest standard error of 6.4*e*–3*

| | | Number of iterations | | | | |
|---|---|---|---|---|---|---|
| **Method** | $b_n =$ | **1e3** | **5e3** | **1e4** | **5e4** | **1e5** |
| BM | $\lfloor n^{1/3} \rfloor$ | 0.9315 | 0.939 | 0.937 | 0.942 | 0.943 |
| Brt | | 0.93 | 0.9395 | 0.936 | 0.9415 | 0.944 |
| OBM | | 0.9305 | 0.9395 | 0.936 | 0.942 | 0.944 |
| TH | | 0.936 | 0.9465 | 0.9395 | 0.9465 | 0.947 |
| BM | $\lfloor n^{1/2} \rfloor$ | 0.9415 | 0.948 | 0.939 | 0.947 | 0.949 |
| Brt | | 0.933 | 0.946 | 0.935 | 0.947 | 0.9475 |
| OBM | | 0.9385 | 0.947 | 0.9355 | 0.9475 | 0.9475 |
| TH | | 0.9365 | 0.9465 | 0.9365 | 0.948 | 0.948 |
| BM | $\lfloor n^{2/3} \rfloor$ | 0.9475 | 0.9445 | 0.9385 | 0.95 | 0.9465 |
| Brt | | 0.9105 | 0.9265 | 0.9275 | 0.9445 | 0.9425 |
| OBM | | 0.9245 | 0.935 | 0.932 | 0.947 | 0.944 |
| TH | | 0.9115 | 0.927 | 0.927 | 0.9435 | 0.9425 |

Comparing the results in Tables 1 and 2 we see that for highly correlated chains, a larger simulation effort and $b_n$ are required to achieve good coverage. However, with a lower correlation, larger values for $b_n$ can result in worse coverage probabilities, especially for smaller simulation efforts.

4.2. *Bayesian probit regression.* Consider the Lupus Data from van Dyk and Meng [51]. This example is concerned with the occurrence of latent membranous lupus nephritis using $y_i$, an indicator of the disease (1 for present), $x_{i1}$, the difference between IgG3 and IgG4 (immunoglobulin G) and $x_{i2}$, IgA (immunoglobulin A) where $i = 1, \dots, 55$. Suppose

$$\Pr(Y_i = 1) = \Phi(\beta_0 + \beta_1 x_{i1} + \beta_2 x_{i2})$$

and assign a flat prior (three-dimensional Lebesgue measure) on $\beta := (\beta_0, \beta_1, \beta_2)$. Roy and Hobert [47] verify that the resulting posterior distribution is proper. Our goal is to estimate the posterior expectation of $\beta$, $E_\pi \beta$.

We will sample from $\pi(\beta|y)$ using the PX-DA algorithm of Liu and Wu [31]. Let $X$ be the $55 \times 3$ design matrix whose $i$th row is $x_i^T = (1, x_{i1}, x_{i2})$ and let $\text{TN}(\mu, \sigma^2, w)$ denote a normal distribution with mean $\mu$ and variance $\sigma^2$ that is truncated to be positive if $w = 1$ and negative if $w = 0$. Then one iteration $\beta \to \beta'$ requires:

TABLE 2

*Table of coverage probabilities for 2000 replications using the* AR*(1) example with* $\rho = 0.95$*. All calculations were based on the nominal level of 0.95. The standard errors for these numbers are easily calculated as* $\sqrt{\hat{p}(1-\hat{p})/2000}$ *which results in a largest standard error of 0.011*

| | | Number of iterations | | | | |
|---|---|---|---|---|---|---|
| **Method** | $b_n =$ | **1e3** | **5e3** | **1e4** | **5e4** | **1e5** |
| BM | $\lfloor n^{1/3} \rfloor$ | 0.614 | 0.738 | 0.766 | 0.842 | 0.872 |
| Brt | | 0.606 | 0.736 | 0.764 | 0.841 | 0.871 |
| OBM | | 0.61 | 0.736 | 0.764 | 0.842 | 0.872 |
| TH | | 0.61 | 0.74 | 0.77 | 0.854 | 0.886 |
| BM | $\lfloor n^{1/2} \rfloor$ | 0.838 | 0.903 | 0.9155 | 0.94 | 0.9425 |
| Brt | | 0.807 | 0.893 | 0.911 | 0.9365 | 0.9385 |
| OBM | | 0.821 | 0.895 | 0.913 | 0.937 | 0.9395 |
| TH | | 0.822 | 0.9055 | 0.9235 | 0.943 | 0.945 |
| BM | $\lfloor n^{2/3} \rfloor$ | 0.927 | 0.9385 | 0.933 | 0.948 | 0.9465 |
| Brt | | 0.872 | 0.916 | 0.9185 | 0.944 | 0.942 |
| OBM | | 0.89 | 0.925 | 0.924 | 0.9455 | 0.943 |
| TH | | 0.885 | 0.92 | 0.924 | 0.9435 | 0.9425 |

1. Draw $z_1, \ldots, z_{55}$ independently with $z_i \sim \mathrm{TN}(x_i^T \beta, 1, y_i)$.
2. Draw $g^2 \sim \mathrm{Gamma}(\frac{55}{2}, \frac{1}{2}\sum_{i=1}^{55}[z_i - x_i^T(X^TX)^{-1}X^Tz]^2)$ and set $z' = (gz_1, \ldots, gz_{55})^T$.
3. Draw $\beta' \sim \mathrm{N}((X^TX)^{-1}X^Tz', (X^TX)^{-1})$.

Roy and Hobert [47] prove that this sampler is geometrically ergodic.

Just as with the AR(1) example, we want to compare observed coverage probabilities based on different methods of estimating the variance of the asymptotic normal distribution. Of course, we need the true value of the unknown posterior expectation $E_\pi \beta$ to do this. To solve this problem, we calculated an estimate of $E_\pi \beta$ from a long simulation (1e8 iterations) of the PX-DA chain and declared the observed sample averages to be the truth. Table 3 shows the observed sample averages along with MCSEs which were calculated using BM with a batch size of $b_n = n^{1/2}$.

We now turn our attention to comparing coverage probabilities in the context of fixed-width methodology, described in Section 1, for each combination of BM, OBM, Brt and TH with sampling plans, $b_n = \lfloor n^\nu \rfloor$ where $\nu \in \{1/3, 1/2\}$. Suppose we want to estimate each component of $E_\pi \beta$ to within $\epsilon \in \{0.1, 0.2, 0.3\}$ while requiring a minimum simulation effort of $n^* \in \{5e4, 1e4, 5e3\}$, respectively. (Recall that the asymptotics require the desired half-width $\epsilon \to 0$ to achieve an asymptotically valid interval.) That is, using the PX-DA algorithm started from the maximum likelihood esti-mate $\hat{\beta} = (-1.778, 4.374, 2.482)$, we terminate the simulation the first time the following inequality holds:

$$\max\left\{t_* \frac{\hat{\sigma}_{\beta_j}}{\sqrt{n}}, \text{ for } j = 0, 1, 2\right\} + \epsilon I(n \le n^*) + n^{-1} \le \epsilon, \tag{12}$$

where $t_*$ is the appropriate critical value for a nominal 95% interval and $\hat{\sigma}_{\beta_j}$ is an estimate for $\sigma_{\beta_j}$ obtained under the settings described above. If (12) was not satisfied, then an additional 10% of the current number of iterations were simulated before checking the criterion again. To estimate the coverage probabilities we obtained 1000 independent replications of this procedure. Table 4 shows the estimated coverage probabilities and the mean number of iterations at termination.

Table 3
*Values treated as the "truth" for estimating confidence interval coverage probabilities based on* 1e8 *iterations*

| | $\boldsymbol{\beta_0}$ | $\boldsymbol{\beta_1}$ | $\boldsymbol{\beta_2}$ |
|---|---|---|---|
| $\hat{\beta}_i$ | −3.0166 | 6.9107 | 3.9792 |
| $\hat{\sigma}_{\beta_i}$ | 11.85 | 22.60 | 14.74 |
| MCSE | 1.18e–3 | 2.26e–3 | 1.47e–3 |

When $\nu = 1/3$, the results are terrible with any method for estimating $\sigma_\beta^2$ or value of $\epsilon$. We can see that even when the simulations are run for more than 1e5 iterations (when $\epsilon = 0.1$), the resulting coverage probabilities are poor. When $\nu = 1/2$, all of the methods result in coverage probabilities slightly lower than the nominal 0.95 level. It also appears that using TH results in slightly better coverage probabilities while requiring slightly more simulation effort. Somewhat surprisingly, the observed coverage probabilities did not uniformly improve as $\epsilon$ decreased.

To this point we have examined the performance of multiple confidence intervals individually but there is an obvious inherent multiplicity issue. Thus, we use a Bonferroni correction to calculate simultaneous confidence intervals. We maintain all settings described above with $\epsilon = 0.2$, except that instead of using nominal 95% we use nominal 98 1/3% confidence intervals. Table 5 shows the estimated coverage probabilities for $E_\pi \beta$ and mean number of iterations at termination based on 1000 independent replications. With $\nu = 1/3$, the results have improved but are still quite poor. However, in the $\nu = 1/2$ case, all individual confidence intervals perform well with observed coverage probabilities close to the nominal 0.9833 level. In addition, the simultaneous intervals have observed coverage probabilities greater than the 0.95 nominal level.

TABLE 4

*Summary of results for using fixed-width methods for the Lupus data Bayesian probit regression. Coverage probabilities using calculated half-width have MCSEs between 1.1e–2 and 1.5e–2 when $b_n = \lfloor n^{1/3} \rfloor$ and between 7.3e–3 and 8.7e–3 when $b_n = \lfloor n^{1/2} \rfloor$. The table also shows the mean (s.e.) simulation effort at termination in terms of number of iterations*

| | | $b_n = \lfloor n^{1/3} \rfloor$ | | | | $b_n = \lfloor n^{1/2} \rfloor$ | | | |
|---|---|---|---|---|---|---|---|---|---|
| **Method** | $\epsilon$ | $\beta_0$ | $\beta_1$ | $\beta_2$ | $n$ | $\beta_0$ | $\beta_1$ | $\beta_2$ | $n$ |
| BM | 0.3 | 0.699 | 0.709 | 0.704 | 6.43e3 (33) | 0.921 | 0.925 | 0.927 | 1.85e4 (89) |
| Brt | | 0.699 | 0.708 | 0.705 | 6.34e3 (32) | 0.926 | 0.923 | 0.930 | 1.81e4 (81) |
| OBM | | 0.700 | 0.710 | 0.706 | 6.37e3 (32) | 0.926 | 0.922 | 0.932 | 1.83e4 (81) |
| TH | | 0.699 | 0.710 | 0.703 | 6.53e3 (34) | 0.936 | 0.938 | 0.941 | 1.95e4 (85) |
| BM | 0.2 | 0.794 | 0.781 | 0.782 | 1.91e4 (60) | 0.922 | 0.929 | 0.930 | 4.54e4 (154) |
| Brt | | 0.797 | 0.775 | 0.781 | 1.89e4 (60) | 0.927 | 0.928 | 0.934 | 4.48e4 (141) |
| OBM | | 0.796 | 0.777 | 0.781 | 1.90e4 (60) | 0.928 | 0.928 | 0.936 | 4.50e4 (140) |
| TH | | 0.804 | 0.790 | 0.794 | 1.97e4 (61) | 0.940 | 0.944 | 0.944 | 4.79e4 (144) |
| BM | 0.1 | 0.854 | 0.849 | 0.853 | 1.11e5 (187) | 0.923 | 0.920 | 0.917 | 1.94e5 (454) |
| Brt | | 0.851 | 0.850 | 0.850 | 1.11e5 (185) | 0.925 | 0.922 | 0.917 | 1.93e5 (410) |
| OBM | | 0.853 | 0.850 | 0.850 | 1.11e5 (185) | 0.925 | 0.922 | 0.917 | 1.93e5 (412) |
| TH | | 0.860 | 0.858 | 0.850 | 1.18e5 (195) | 0.932 | 0.929 | 0.932 | 2.02e5 (419) |

4.2.1. *Comparison to regeneration.* Finally, we compare BM, Brt, OBM and TH with regenerative simulation (RS) in terms of coverage probabilities. RS is often viewed as the gold standard; see the references in Section 1. We are still using the model given above and sampling via the PX-DA

TABLE 5

*Summary of results for using fixed-width methods with a Bonferroni correction for the Lupus data Bayesian probit regression. Coverage probabilities using calculated half-width have MCSEs of between 9e–3 and 1e–2 when $b_n = \lfloor n^{1/3} \rfloor$ and between 4.9e–3 and 6e–3 when $b_n = \lfloor n^{1/2} \rfloor$. The table also shows the mean (s.e.) simulation effort at termination in terms of number of iterations*

| **Method** | $b_n =$ | $\beta_0$ | $\beta_1$ | $\beta_2$ | **Simultaneous** | $n$ |
|---|---|---|---|---|---|---|
| BM | $\lfloor n^{1/3} \rfloor$ | 0.911 | 0.904 | 0.906 | 0.878 | 3.21e4 (85) |
| Brt | | 0.909 | 0.904 | 0.909 | 0.879 | 3.19e4 (85) |
| OBM | | 0.910 | 0.905 | 0.909 | 0.879 | 3.20e4 (85) |
| TH | | 0.917 | 0.909 | 0.911 | 0.884 | 3.34e4 (87) |
| BM | $\lfloor n^{1/2} \rfloor$ | 0.973 | 0.975 | 0.972 | 0.965 | 6.97e4 (210) |
| Brt | | 0.973 | 0.971 | 0.974 | 0.965 | 6.86e4 (189) |
| OBM | | 0.973 | 0.970 | 0.972 | 0.963 | 6.90e4 (191) |
| TH | | 0.973 | 0.976 | 0.972 | 0.969 | 7.30e4 (195) |

TABLE 6
*Coverage probabilities comparing BM, OBM, and SV using 7e5 iterations and $b_n = \lfloor n^{1/2} \rfloor$ to RS. MCSEs vary between 6.9e–3 and 7.6e–3*

| **Method** | $\boldsymbol{\beta_0}$ | $\boldsymbol{\beta_1}$ | $\boldsymbol{\beta_2}$ |
|---|---|---|---|
| BM | 0.945 | 0.945 | 0.950 |
| Brt | 0.941 | 0.942 | 0.948 |
| OBM | 0.942 | 0.942 | 0.948 |
| TH | 0.946 | 0.945 | 0.950 |
| RS | 0.938 | 0.938 | 0.938 |

algorithm, but we are not concerned with fixed-width methodology in this section. Roy and Hobert [47] implement RS for this example and we use their settings exactly except that we use only 50 regenerations. We obtained 1000 independent replications resulting in a mean simulation effort of 7.12e5 (3.2e3). For a fair comparison in terms of simulation effort, confidence intervals for $E_\pi \beta$ were calculated using BM, OBM and SV using $b_n = \lfloor n^{1/2} \rfloor$ from 1000 independent replications with 7e5 iterations each. Table 6 shows the resulting coverage probabilities for $E_\pi \beta$. The coverage probabilities from BM, Brt, OBM and TH are all slightly larger than those of RS; however, the results from all the methods are within two standard errors of the nominal 0.95 level.

4.3. *Summary.* In our examples, we consider different estimators of $\sigma_g^2$. Roughly all of the methods considered resulted in similar performance in terms of estimated coverage probabilities. Recall, OBM and Brt are asymptotically equivalent and the simulation results show there is little difference between the two in finite samples. However, the TH estimator tends to perform slightly better than OBM and Brt. In our experience, Brt (and SV methods in general) tended to be slightly faster than OBM from a computational perspective. As suggested by the theoretical results in Section 3.3, the estimator from BM was observed to be more variable than OBM. This was consistent in both examples in multiple realizations of the simulation.

Using the Bayesian probit regression model we compared our methods to RS. The resulting simulation showed all of the methods performed very well. The advantage of RS is that in the fixed-width setting the actual chain does not need to be stored as the simulation progresses. However, RS requires a theoretical cost that, while not overly burdensome in our view, may dissuade the typical practitioner. The resulting simulation is also dependent on the length of the regeneration tours which can be extremely long in even moderately large finite state spaces or as the dimension of the Markov chain increases (Gilks, Roberts and Sahu [14]) or in variable-at-a-time Metropolis–Hastings implementations (Neath and Jones [39]). In contrast, BM, OBM and SV are relatively simple to implement though they can require saving the entire chain if fixed-width methods are employed. Given the current price of computer memory, this is clearly not the obstacle it was in the past.

Another simulation goal was to investigate the finite sample behavior of different batch size selection. Theoretically, we showed that the batch size should increase at a rate proportional to $n^{1/3}$ but the proportionality constant is unknown. In our examples, using $b_n = \lfloor n^{1/3} \rfloor$ seemed to give very poor results because the batch size or truncation point was too small. In realistic examples with higher correlations, the larger batch size $b_n = \lfloor n^{1/2} \rfloor$ worked well agreeing with the previous work of Jones et al. [24]. Our investigation of $b_n = \lfloor n^{2/3} \rfloor$ worked well in high correlation settings, though for long chains more computational effort was required.

On balance, we would recommend using the SV method with the Tukey–Hanning window with $b_n = \lfloor n^{1/2} \rfloor$ as a default method. If the required moment conditions for the SV methods were too much, then we would employ OBM.

## APPENDIX A: BROWNIAN MOTION

Let $B = \{B(t), t \geq 0\}$ denote a standard Brownian motion. Define $\bar{B}_j(k) := k^{-1}(B(j+k) - B(j))$ and $\bar{B}_n := n^{-1}B(n)$. We will require the following two results from Csörgő and Révész [8] on the increments of Brownian motion.

LEMMA 1. *For all $\epsilon > 0$ and for almost all sample paths there exists $n_0(\epsilon)$ such that for all $n \geq n_0$,*

$$|B(n)| < (1+\epsilon)[2n \log\log n]^{1/2}.$$

LEMMA 2. *Suppose Assumption 2 holds, then for all $\epsilon > 0$ and for almost all sample paths, there exists $n_0(\epsilon)$ such that for all $n \geq n_0$,*

$$\sup_{0 \leq t \leq n-b_n} \sup_{0 \leq s \leq b_n} |B(t+s) - B(t)| < (1+\epsilon)\left[2b_n\left(\log \frac{n}{b_n} + \log\log n\right)\right]^{1/2}.$$

A *strong invariance principle* holds if there exists a nonnegative increasing function $\psi(n)$ on the positive integers, a constant $0 < \sigma_g < \infty$ and a sufficiently rich probability space $\Omega$ such that

$$\left|\sum_{i=1}^{n} g(X_i) - nE_\pi g - \sigma_g B(n)\right| = O(\psi(n)) \qquad \text{w.p. 1 as } n \to \infty, \tag{13}$$

where the w.p. 1 in (13) means for almost all sample paths. Alternatively, (13) can be expressed as there exists $n_0$ and a finite random variable $C$ such that for almost all $\omega \in \Omega$,

$$\left|\sum_{i=1}^{n} Y_i - nE_\pi g - \sigma_g B(n)\right| < C(\omega)\psi(n) \tag{14}$$

for all $n > n_0$. A strong invariance principle is enough to guarantee both a strong law (1), a central limit theorem (2) and a functional central limit theorem among other properties (see Philipp and Stout [40] and Damerdji [9]). We will rely on the following result to connect the strong invariance principle to the convergence rate of a Harris ergodic Markov chain (see Bednorz and Latuszyński [2] and Jones et al. [24] for a proof).

LEMMA 3. *Let $g : \mathsf{X} \mapsto \mathbb{R}$ be a Borel function and let $X$ be a geometrically ergodic Markov chain with invariant distribution $\pi$. Suppose Assumption 3 holds and $E_\pi|g|^{2+\delta+\epsilon} < \infty$ for some $\delta > 0$ and some $\epsilon > 0$, then a strong invariance principle holds with $\psi(n) = n^\alpha \log n$ where $\alpha = 1/(2+\delta)$.*

## APPENDIX B: STRONG CONSISTENCY PROOFS

**B.1. Proof of Theorem 1.** The proof will be constructed in 3 stages given in lemmas below, but first we define some notation. Recall that $X = \{X_1, X_2, \ldots\}$ is a Harris ergodic Markov chain with invariant distribution $\pi$ and $Y_i = g(X_i) - E_\pi g$ for $i = 1, 2, 3, \ldots$. Further define $\bar{Y}_j(k) = k^{-1}\sum_{i=1}^{k} Y_{j+i}$ for $j = 0, \ldots, n-b_n$ and $k = 1, \ldots, b_n$ and $\bar{Y}_n = n^{-1}\sum_{i=1}^{n} Y_i$. Next let

$$\hat{\sigma}^2_{w,n} := \frac{1}{n}\sum_{j=0}^{n-b_n}\sum_{k=1}^{b_n} k^2 \Delta_2 w_n(k)[\bar{Y}_j(k) - \bar{Y}_n]^2$$

and

$$\tilde{\sigma}^2_* := \frac{1}{n}\sum_{j=0}^{n-b_n}\sum_{k=1}^{b_n} k^2 \Delta_2 w_n(k)[\bar{B}_j(k) - \bar{B}_n]^2.$$

LEMMA 4. *Suppose (13) holds with $\psi(n) = n^\alpha \log n$ where $\alpha = 1/(2+\delta)$ and Assumptions 1 and 2 hold. If, as $n \to \infty$,*

$$b_n^{1/2} n^\alpha (\log n)^{3/2} \sum_{k=1}^{b_n} |\Delta_2 w_n(k)| \to 0 \quad \textit{and} \tag{15}$$

$$n^{2\alpha}(\log n)^2 \sum_{k=1}^{b_n} |\Delta_2 w_n(k)| \to 0, \tag{16}$$

*then* $\hat{\sigma}^2_{w,n} - \sigma_g^2 \tilde{\sigma}^2_* \to 0$ *w.p. 1.*

PROOF. Notice

$$\hat{\sigma}^2_{w,n} - \sigma_g^2 \tilde{\sigma}^2_* = \frac{1}{n} \sum_{j=0}^{n-b_n} \sum_{k=1}^{b_n} k^2 \Delta_2 w_n(k)([\bar{Y}_j(k) - \bar{Y}_n]^2 - \sigma_g^2 [\bar{B}_j(k) - \bar{B}_n]^2).$$

Let $A_k = k[\bar{Y}_j(k) - \sigma_g \bar{B}_j(k)]$, $D_k = B(j+k) - B(j)$, $E_{n,k} = k\bar{B}_n$ and $F_{n,k} = k[\bar{Y}_n - \sigma_g \bar{B}_n]$. Then

$$\begin{aligned} k[\bar{Y}_j(k) - \bar{Y}_n] &= k[\bar{Y}_j(k) - \bar{Y}_n \pm \sigma_g \bar{B}_j(k) \pm \sigma_g \bar{B}_n] \\ &= k[\bar{Y}_j(k) - \sigma_g \bar{B}_j(k)] + \sigma_g k \bar{B}_j(k) \\ &\quad - \sigma_g k \bar{B}_n - k[\bar{Y}_n - \sigma_g \bar{B}_n] \\ &= k[\bar{Y}_j(k) - \sigma_g \bar{B}_j(k)] + \sigma_g [B(j+k) - B(j)] \\ &\quad - \sigma_g k \bar{B}_n - k[\bar{Y}_n - \sigma_g \bar{B}_n] \\ &= A_k + \sigma_g (D_k - E_{n,k}) - F_{n,k}. \end{aligned}$$

Hence

$$\begin{aligned} &|\hat{\sigma}^2_{w,n} - \sigma_g^2 \tilde{\sigma}^2_*| \\ &\quad \le \frac{1}{n} \sum_{j=0}^{n-b_n} \sum_{k=1}^{b_n} |\Delta_2 w_n(k) [(A_k + \sigma_g (D_k - E_{n,k}) - F_{n,k})^2 \\ &\qquad\qquad - \sigma_g^2 (D_k - E_{n,k})^2]| \qquad (17) \\ &\quad \le \frac{1}{n} \sum_{j=0}^{n-b_n} \sum_{k=1}^{b_n} |\Delta_2 w_n(k)| (A_k^2 + F_{n,k}^2 + 2\sigma_g |A_k D_k| + 2\sigma_g |A_k E_{n,k}| \\ &\qquad\qquad + 2|A_k F_{n,k}| + 2\sigma_g |D_k F_{n,k}| + 2\sigma_g |E_{n,k} F_{n,k}|). \end{aligned}$$

It suffices to show that each of the 7 sums in (17) tend to 0 as $n \to \infty$. By our assumption of (13) we have that for sufficiently large $n$,

$$\left| \sum_{i=1}^{n} Y_i - \sigma_g B(n) \right| \le C n^{\alpha} \log n. \qquad (18)$$

1. From (18), we obtain

$$\begin{aligned} |A_k| &= \left| \left( \sum_{i=1}^{j+k} Y_i - \sigma_g B(j+k) \right) - \left( \sum_{i=1}^{j} Y_i - \sigma_g B(j) \right) \right| \\ &\le C(j+k)^{\alpha} \log(j+k) + C(j)^{\alpha} \log(j), \end{aligned}$$

and since $j + k \leq n$,

$$|A_k| \leq 2Cn^{\alpha} \log n. \tag{19}$$

Hence

$$\frac{1}{n}\sum_{j=0}^{n-b_n}\sum_{k=1}^{b_n}|\Delta_2 w_n(k)|A_k^2 \leq 4C^2 n^{2\alpha}(\log n)^2 \sum_{k=1}^{b_n}|\Delta_2 w_n(k)| \to 0$$

as $n \to \infty$ by (16).

2. From (18) and the fact that $k \leq b_n \leq n$

$$|F_{n,k}| \leq Ckn^{\alpha-1}\log n, \tag{20}$$

resulting in

$$\begin{aligned}\frac{1}{n}\sum_{j=0}^{n-b_n}\sum_{k=1}^{b_n}|\Delta_2 w_n(k)|F_{n,k}^2 &\leq C^2 n^{2\alpha-2}(\log n)^2 \sum_{k=1}^{b_n} k^2|\Delta_2 w_n(k)| \\ &\leq C^2 b_n^2 n^{2\alpha-2}(\log n)^2 \sum_{k=1}^{b_n}|\Delta_2 w_n(k)| \to 0\end{aligned}$$

as $n \to \infty$ by (16).

3. From Lemma 2,

$$\begin{aligned}|D_k| &= |B(j+k) - B(j)| \\ &\leq \sup_{0\leq t\leq n-b_n}\ \sup_{0\leq s\leq b_n}|B(t+s) - B(t)| \\ &\leq (1+\epsilon)\left(2b_n\left(\log\frac{n}{b_n} + \log\log n\right)\right)^{1/2} \\ &\leq 2(1+\epsilon)b_n^{1/2}(\log n)^{1/2}.\end{aligned} \tag{21}$$

Combining (19) with (21), we obtain

$$\begin{aligned}&\frac{1}{n}\sum_{j=0}^{n-b_n}\sum_{k=1}^{b_n}|\Delta_2 w_n(k)|2\sigma_g|A_k D_k| \\ &\quad \leq 8C\sigma_g(1+\epsilon)b_n^{1/2}n^{\alpha}(\log n)^{3/2}\sum_{k=1}^{b_n}|\Delta_2 w_n(k)| \to 0\end{aligned}$$

as $n \to \infty$ by (15).

4. From Lemma 1,

$$|E_{n,k}| \leq \sqrt{2}(1+\epsilon)kn^{-1/2}(\log\log n)^{1/2}. \tag{22}$$

Combining (19) with (22), get

$$\frac{1}{n}\sum_{j=0}^{n-b_n}\sum_{k=1}^{b_n}|\Delta_2 w_n(k)|2\sigma_g|A_k E_{n,k}|$$

$$\le 2^{5/2}C\sigma_g(1+\epsilon)n^{\alpha-1/2}\log n(\log\log n)^{1/2}\sum_{k=1}^{b_n}k|\Delta_2 w_n(k)|$$

$$\le 2^{5/2}C\sigma_g(1+\epsilon)(b_n/n)^{1/2}b_n^{1/2}n^{\alpha}(\log n)^{3/2}\sum_{k=1}^{b_n}|\Delta_2 w_n(k)|\to 0$$

as $n\to\infty$ by (15) and Assumption 2.

5. From (19) and (20) we have

$$\frac{1}{n}\sum_{j=0}^{n-b_n}\sum_{k=1}^{b_n}|\Delta_2 w_n(k)|2|A_k F_{n,k}|\le 4C^2 n^{2\alpha-1}(\log n)^2\sum_{k=1}^{b_n}k|\Delta_2 w_n(k)|\to 0$$

as $n\to\infty$ by (16) and Assumption 2.

6. From (20) and (21),

$$\frac{1}{n}\sum_{j=0}^{n-b_n}\sum_{k=1}^{b_n}|\Delta_2 w_n(k)|2\sigma_g|D_k F_{n,k}|$$

$$\le 4\sigma_g(1+\epsilon)b_n^{1/2}n^{\alpha-1}(\log n)^{3/2}\sum_{k=1}^{b_n}k|\Delta_2 w_n(k)|\to 0$$

as $n\to\infty$ by (15) and Assumption 2.

7. From (20) and (22),

$$\frac{1}{n}\sum_{j=0}^{n-b_n}\sum_{k=1}^{b_n}|\Delta_2 w_n(k)|2\sigma_g|E_{n,k}F_{n,k}|$$

$$\le 2^{3/2}C\sigma_g(1+\epsilon)n^{\alpha-3/2}\log n(\log\log n)^{1/2}\sum_{k=1}^{b_n}k^2|\Delta_2 w_n(k)|\to 0$$

as $n\to\infty$ by (15) and Assumption 2. □

LEMMA 5. *Let $X$ be a geometrically ergodic Markov chain with invariant distribution $\pi$ and $g:\mathsf{X}\to\mathbb{R}$ be a Borel function with $E_\pi|g|^{4+\delta+\epsilon}<\infty$ for some $\delta>0$ and $\epsilon>0$. Set $h(X_i)=[g(X_i)-E_\pi g]^2$ for $i\ge 1$. If Assumptions 2 and 3 hold and $b_n^{-1}n^{2\alpha}\log n\to 0$ as $n\to\infty$ where $\alpha=1/(4+\delta)$, then $b_n^{-1}\sum_{i=1}^{b_n}h(X_i)$ and $b_n^{-1}\sum_{i=n-b_n+1}^{n}h(X_i)$ stay bounded as $n\to\infty$ w.p. 1.*

PROOF. The ergodic theorem implies that, w.p. 1, $b_n^{-1}\sum_{i=1}^{b_n} h(X_i)$ converges to a finite limit and hence stays bounded as $n \to \infty$ w.p. 1.

Note that if $E_\pi|g|^{4+\delta+\epsilon} < \infty$, then $E_\pi|h|^{2+\delta/2+\epsilon/2} < \infty$, and hence our assumptions with Lemma 3 yield the existence of $n_0$ such that if $n > n_0$,

$$\left|\sum_{i=1}^{n} h(X_i) - nE_\pi h(X_1) - \sigma_h B(n)\right| < Cn^{2\alpha}\log n, \tag{23}$$

where $\alpha = 1/(4+\delta)$.

Next, for all $\varepsilon > 0$ and sufficiently large $n(\varepsilon)$,

$$\begin{aligned}
&\frac{1}{b_n}\left|\sum_{i=n-b_n+1}^{n} h(X_i)\right| \\
&\quad = \frac{1}{b_n}\left|\sum_{i=1}^{n} h(X_i) - \sum_{i=1}^{n-b_n} h(X_i)\right| \\
&\quad = \frac{1}{b_n}\left|\left(\sum_{i=1}^{n} h(X_i) - nE_\pi h(X_1) - \sigma_h B(n)\right)\right. \\
&\qquad\qquad - \left(\sum_{i=1}^{n-b_n} h(X_i) - (n-b_n)E_\pi h(X_1) - \sigma_h B(n-b_n)\right) \\
&\qquad\qquad\qquad \left. + \sigma_h(B(n) - B(n-b_n)) + b_n E_\pi h(X_1)\right| \\
&\quad \le \frac{1}{b_n}\left[\left|\sum_{i=1}^{n} h(X_i) - nE_\pi h(X_1) - \sigma_h B(n)\right|\right. \\
&\qquad\qquad + \left|\sum_{i=1}^{n-b_n} h(X_i) - (n-b_n)E_\pi h(X_1) - \sigma_h B(n-b_n)\right| \\
&\qquad\qquad\qquad \left. + \sigma_h|B(n) - B(n-b_n)| + b_n E_\pi h(X_1)\right] \\
&\quad \le \frac{1}{b_n}\left[2Cn^{2\alpha}\log n + (1+\varepsilon)\left(2b_n\left(\log\frac{n}{b_n} + \log\log n\right)\right)^{1/2}\right] + E_\pi h(X_1) \\
&\quad = E_\pi h(X_1) + 2Cb_n^{-1}n^{2\alpha}\log n + O((b_n^{-1}\log n)^{1/2}),
\end{aligned}$$

where the second inequality follows from (23) and Lemma 2. Hence, $b_n^{-1} \times |\sum_{i=n-b_n+1}^{n} h(X_i)|$ stays bounded w.p. 1 since $b_n^{-1}n^{2\alpha}\log n \to 0$ as $n \to \infty$. □

LEMMA 6. *Suppose Assumptions 1 and 2 hold. Further assume the conditions of Lemma 5 and that:*

1. *there exists a constant* $c \geq 1$ *such that* $\sum_n (b_n/n)^c < \infty$;

2. $$b_n n^{-1} \sum_{k=1}^{b_n} k|\Delta_1 w_n(k)| \to 0 \qquad \text{as } n \to \infty \quad \text{and}$$

3. $b_n n^{-1} \log n \to 0$ *as* $n \to \infty$.

*Then there exists a sequence of random variables* $d_n$ *such that*

$$\hat{\sigma}^2_{w,n} + d_n = \hat{\sigma}^2_S,$$

*and* $d_n \to 0$ *w.p. 1 as* $n \to \infty$. *Also,* $\tilde{\sigma}^2_* \to 1$ *as* $n \to \infty$.

PROOF. This follows immediately from the conclusion of Lemma 5 and results in Damerdji [9], page 1430. □

PROOF OF THEOREM 1. The result follows by combining Lemmas 3, 4, and 6. Note that since $X$ is geometrically ergodic and $E_\pi|g|^{4+\delta+\epsilon} < \infty$ for some $\delta > 0$ and $\epsilon > 0$, the conclusion of Lemma 4 holds with $\alpha = 1/(4+\delta)$. □

**B.2. A condition for lag windows.**

LEMMA 7. *Suppose* $w$ *is defined on* $[0,1]$ *such that* $w(0) = 1$ *and* $w(1) = 0$. *Further assume that* $w$ *is twice continuously differentiable. Also, assume that* $D_1$ *and* $D_2$ *are finite constants such that* $|w'(x)| \leq D_1$ *and* $|w''(x)| \leq D_2$. *Then as* $n \to \infty$, *condition* (a) *of Theorem 1 holds if* $b_n^2 n^{-1} \to 0$ *while condition* (d) *of Theorem 1 is satisfied if* $b_n^{-1} n^{2\alpha} (\log n)^3 \to 0$.

PROOF. Suppose $1 \leq k \leq b_n - 1$ and let

$$\Delta_2 w\left(\frac{k}{b_n}\right) = w\left(\frac{k-1}{b_n}\right) - 2w\left(\frac{k}{b_n}\right) + w\left(\frac{k+1}{b_n}\right)$$
$$= \left[w\left(\frac{k-1}{b_n}\right) - w\left(\frac{k}{b_n}\right)\right] - \left[w\left(\frac{k}{b_n}\right) - w\left(\frac{k+1}{b_n}\right)\right].$$

The mean value theorem guarantees the existence of $c_1 \in (\frac{k-1}{b_n}, \frac{k}{b_n})$ and $c_2 \in (\frac{k}{b_n}, \frac{k+1}{b_n})$ so that

$$\Delta_2 w\left(\frac{k}{b_n}\right) = -\frac{w'(c_1)}{b_n} + \frac{w'(c_2)}{b_n}.$$

A second application of the mean value theorem yields a constant $c \in (c_1, c_2)$ such that

$$\Delta_2 w\left(\frac{k}{b_n}\right) = \frac{w''(c)}{b_n}(c_2 - c_1).$$

Since $c_2 - c_1 \leq 2/b_n$ and $|w''(x)| \leq D_2$ we have

$$\left|\Delta_2 w\left(\frac{k}{b_n}\right)\right| \leq \frac{2D_2}{b_n^2}.$$

Then

$$\begin{aligned}
\sum_{k=1}^{b_n}\left|\Delta_2 w\left(\frac{k}{b_n}\right)\right| &= \sum_{k=1}^{b_n-1}\left|\Delta_2 w\left(\frac{k}{b_n}\right)\right| + \left|w\left(\frac{b_n-1}{b_n}\right) - w(1)\right| \\
&\leq \sum_{k=1}^{b_n-1}\frac{2D_2}{b_n^2} + \frac{D_1}{b_n} \\
&= \frac{2(b_n-1)D_2}{b_n^2} + \frac{D_1}{b_n}.
\end{aligned}$$

Next observe that

$$\begin{aligned}
b_n n^{2\alpha}(\log n)^3\left(\sum_{k=1}^{b_n}\left|\Delta_2 w\left(\frac{k}{b}\right)\right|\right)^2 &\leq b_n n^{2\alpha}(\log n)^3\left[\frac{2(b_n-1)D_2}{b_n^2} + \frac{D_1}{b_n}\right]^2 \\
&= b_n^{-1} n^{2\alpha}(\log n)^3\left[\frac{2(b_n-1)D_2}{b_n} + D_1\right]^2 \\
&\to 0 \qquad \text{as } n \to \infty,
\end{aligned}$$

since $b_n^{-1} n^{2\alpha}(\log n)^3$ as $n \to \infty$ by assumption. The proofs that our conditions also imply the remaining portion of condition (d) and condition (a) are similar and hence are omitted. □

**B.3. Proof of Theorem 2.** Consider the modified Bartlett lag window defined in Remark 3. Then $\Delta_2 w_n(b_n) = b_n^{-1}$ and $\Delta_2 w_n(k) = 0$ for $k = 1, 2, \dots, b_n - 1$ so that

$$\hat{\sigma}_{w,n}^2 = n^{-1}\sum_{j=0}^{n-b_n}\sum_{k=1}^{b_n} k^2 \Delta_2 w_n(k)[\bar{Y}_j(k) - \bar{Y}_n]^2 = \frac{b_n}{n}\sum_{j=0}^{n-b_n}[\bar{Y}_j(b_n) - \bar{Y}_n]^2 \tag{24}$$

and

$$\tilde{\sigma}_*^2 = \frac{b_n}{n}\sum_{j=0}^{n-b_n}[\bar{B}_j(b_n) - \bar{B}_n]^2. \tag{25}$$

Further, (15) and (16) are satisfied if $b_n^{-1}n^{2\alpha}(\log n)^3 \to 0$ as $n \to \infty$ where $\alpha = 1/(2+\delta)$, and hence our assumptions imply the conclusion of Lemma 4. Since

$$\hat{\sigma}^2_{\text{OBM}} = \frac{nb_n}{(n-b_n)(n-b_n+1)} \sum_{j=0}^{n-b_n} [\bar{Y}_j(b_n) - \bar{Y}_n]^2,$$

the conclusion would follow from Lemmas 3 and 4 if $\tilde{\sigma}^2_* \to 1$ w.p. 1 as $n \to \infty$ because $\hat{\sigma}^2_{w,n}$ is asymptotically equivalent to $\hat{\sigma}^2_{\text{OBM}}$: as $n \to \infty$

$$\frac{\hat{\sigma}^2_{w,n}}{\hat{\sigma}^2_{\text{OBM}}} = \left(1 - \frac{b_n}{n}\right)\left(1 - \frac{b_n}{n} + \frac{1}{b_n}\right) \to 1.$$

Define $U_i := B(i) - B(i-1)$, the increments of Brownian motion, so that $U_1, \ldots, U_n$ are i.i.d. $\text{N}(0,1)$. Further define $T_i = U_i - \bar{B}_n$ for $i = 1, \ldots, n$ and set

$$\begin{aligned}\tilde{d}_n = \frac{1}{n}\Bigg[\sum_{l=1}^{b_n} \frac{1}{b_n}\left(\sum_{i=1}^{l-1} T_i^2 + \sum_{i=n-b_n+l+1}^{n} T_i^2\right)\Bigg] \\ + 2\sum_{s=1}^{b_n-1}\Bigg[\sum_{l=1}^{b_n-s} \frac{1}{b_n}\left(\sum_{i=1}^{l-1} T_i T_{s+i} + \sum_{i=n-b_n+l+1}^{n-s} T_i T_{s+i}\right)\Bigg],\end{aligned}$$

where any empty sums are defined to be zero. Letting $\tilde{\gamma}_n(i) = \tilde{\gamma}_n(-i) := n^{-1}\sum_{t=1}^{n-i}(U_t - \bar{B}_n)(U_{t+i} - \bar{B}_n)$ Damerdji [9] shows that under Assumption 2,

$$\tilde{\sigma}^2_* + \tilde{d}_n = \sum_{s=-(b_n-1)}^{b_n-1} \left(1 - \frac{|k|}{b_n}\right)\tilde{\gamma}_n(s).$$

Moreover, under conditions (a) and (b) of Theorem 2, $\tilde{\sigma}^2_* + \tilde{d}_n \to 1$ as $n \to \infty$ w.p. 1. Thus the proof will be complete if we can show that $\tilde{d}_n \to 0$ as $n \to \infty$ w.p. 1.

LEMMA 8. *Suppose Assumption 2 holds. If:*

1. *there exists a constant* $c \geq 1$ *such that* $\sum_n (b_n/n)^{2c} < \infty$;
2. *there exists an integer* $n_0$ *and a constant* $c_1$ *such that for all* $n \geq n_0$ *we have* $\log n / b_n \leq c_1$ *and*
3. *as* $n \to \infty$ *we have* $b_n^4 n^{-3} \log\log n \to 0$ *and* $b_n^2 n^{-2} \log\log n \to 0$,

*then* $\tilde{d}_n \to 0$ *w.p. 1 as* $n \to \infty$.

Before we present the proof we recall two results which will be used below.

LEMMA 9 (Kendall and Stuart [27]). *If* $Z \sim \chi_v^2$, *then for all positive integers* $r$ *there exists a constant* $K := K(r)$ *such that* $E[(Z-v)^{2r}] \leq K v^r$.

LEMMA 10 (Whittle [53]). *Let* $U_1, \dots, U_n$ *be i.i.d. standard normal variables and* $A = \sum_j \sum_k a_{jk} U_j U_k$ *where* $a_{jk}$ *are real coefficients; then for* $c \geq 1$ *we have*

$$E[|A - EA|^{2c}] \leq K(c) \left( \sum_j \sum_k a_{jk}^2 \right)^c$$

*for some constant* $K_c < \infty$ *depending only on* $c$.

PROOF OF LEMMA 8. Note

$$|\tilde{d}_n| \leq \frac{1}{n} \left| \sum_{l=1}^{b_n} b_n^{-1} \left( \sum_{i=1}^{l-1} T_i^2 + \sum_{i=n-b_n+l+1}^{n} T_i^2 \right) \right| \tag{26}$$

$$+ \frac{1}{n} \left| 2 \sum_{s=1}^{b_n-1} \sum_{l=1}^{b_n-s} b_n^{-1} \left( \sum_{i=1}^{l-1} T_i T_{s+i} + \sum_{i=n-b_n+l+1}^{n-s} T_i T_{s+i} \right) \right|. \tag{27}$$

We will show that (26) and (27) each converge almost surely to zero implying the desired result.

First consider (26):

$$\begin{aligned}
&\frac{1}{n} \left| \sum_{l=1}^{b_n} \frac{1}{b_n} \left( \sum_{i=1}^{l-1} T_i^2 + \sum_{i=n-b_n+l+1}^{n} T_i^2 \right) \right| \\
&\quad \leq \frac{1}{n} \left| \sum_{l=1}^{b_n} \frac{1}{b_n} \left( \sum_{i=1}^{b_n-1} T_i^2 + \sum_{i=n-b_n+2}^{n} T_i^2 \right) \right| \\
&\quad = \frac{1}{n} \left| \sum_{i=1}^{b_n-1} T_i^2 + \sum_{i=n-b_n+2}^{n} T_i^2 \right| \\
&\quad = \frac{1}{n} \left| \sum_{i=1}^{b_n-1} (U_i - \bar{B}_n)^2 + \sum_{i=n-b_n+2}^{n} (U_i - \bar{B}_n)^2 \right| \\
&\quad \leq \frac{1}{n} \left| \sum_{i=1}^{b_n-1} U_i^2 + \sum_{i=n-b_n+2}^{n} U_i^2 \right| \\
&\qquad + \frac{1}{n} \left| 2\bar{B}_n \left( \sum_{i=1}^{b_n-1} U_i + \sum_{i=n-b_n+2}^{n} U_i \right) \right| \\
&\qquad + \frac{1}{n} |2(b_n - 1)\bar{B}_n^2|.
\end{aligned}$$

Now we show that each of the three terms above tend to zero.

1. Since $U_1^2, \dots, U_n^2$ are i.i.d. $\chi_1^2$, $\sum_{i=1}^{b_n-1} U_i^2 + \sum_{i=n-b_n+2}^{n} U_i^2 \sim \chi^2_{2(b_n-1)}$. By Lemma 9 we have for every positive integer $r$ that

$$E\left\{\left[\frac{1}{n}\left(\sum_{i=1}^{b_n-1} U_i^2 + \sum_{i=n-b_n+2}^{n} U_i^2\right) - \frac{2(b_n-1)}{n}\right]^{2r}\right\} \le K\left(\frac{2(b_n-1)}{n^2}\right)^r.$$

Now choose $r > 1$ so that

$$\sum_n K\left(\frac{2(b_n-1)}{n^2}\right)^r \le K\sum_n \left(\frac{2}{n}\right)^r < \infty.$$

Let $\epsilon > 0$ be arbitrary. Then by Markov's inequality we have

$$\epsilon \sum_n \Pr\left\{\left[\frac{1}{n}\left(\sum_{i=1}^{b_n-1} U_i^2 + \sum_{i=n-b_n+2}^{n} U_i^2\right) - \frac{2(b_n-1)}{n}\right]^{2r} \ge \epsilon\right\} < \infty.$$

Now a standard Borel–Cantelli argument (Billingsley [3], pages 59 and 60) yields

$$\left[\frac{1}{n}\left(\sum_{i=1}^{b_n-1} U_i^2 + \sum_{i=n-b_n+2}^{n} U_i^2\right) - \frac{2(b_n-1)}{n}\right]^{2r} \to 0 \qquad \text{w.p. 1 as } n \to \infty.$$

Hence

$$\frac{1}{n}\left|\sum_{i=1}^{b_n-1} U_i^2 + \sum_{i=n-b_n+2}^{n} U_i^2\right| \to 0 \qquad \text{w.p. 1 as } n \to \infty,$$

since $b_n/n \to 0$ as $n \to \infty$.

2. Notice

$$\frac{1}{n}\left|2\bar{B}_n\left(\sum_{i=1}^{b_n-1} U_i + \sum_{i=n-b_n+2}^{n} U_i\right)\right| \le \frac{2|\bar{B}_n|}{n}\sum_{i=1}^{n}|U_i|.$$

Since the $|U_i|$ are i.i.d. following a half-normal distribution, the classical strong law implies $n^{-1}\sum_{i=1}^{n}|U_i| \to \sqrt{2/\pi}$ w.p. 1. (Recall that $E|U_i| = \sqrt{2/\pi}$ and $\operatorname{Var}|U_i| = 1 - 2/\pi$.) Combining this with Lemma 1 yields for every $\epsilon > 0$ and sufficiently large $n$,

$$\frac{2|\bar{B}_n|}{n}\sum_{i=1}^{n}|U_i| \le 2^{3/2}(1+\epsilon)\left[\frac{\log\log n}{n}\right]^{1/2}\frac{1}{n}\sum_{i=1}^{n}|U_i| \to 0$$

with probability 1 as $n \to \infty$.

3. Using Lemma 1 shows that for every $\epsilon > 0$ and large enough $n$,

$$\frac{1}{n}|2(b_n-1)\bar{B}_n^2| \le 4(1+\epsilon)^2 \frac{(b_n-1)\log\log n}{n^2} \to 0$$

since $b_n/n \to 0$ as $n \to \infty$.

This establishes that the term in (26) converges almost surely to 0 as $n \to \infty$. Next, consider (27):

$$\begin{aligned}
&\frac{1}{n}\left|2\sum_{s=1}^{b_n-1}\sum_{l=1}^{b_n-s}\frac{1}{b_n}\left(\sum_{i=1}^{l-1}T_iT_{s+i} + \sum_{i=n-b_n+l+1}^{n-s}T_iT_{s+i}\right)\right| \\
&\quad = \frac{1}{n}\left|2\sum_{s=1}^{b_n-1}\sum_{l=1}^{b_n-s}\frac{1}{b_n}\left(\sum_{i=1}^{l-1}(U_i-\bar{B}_n)(U_{s+i}-\bar{B}_n)\right.\right. \\
&\qquad\qquad \left.\left. + \sum_{i=n-b_n+l+1}^{n-s}(U_i-\bar{B}_n)(U_{s+i}-\bar{B}_n)\right)\right| \\
&\quad \le \frac{1}{n}\left|2\sum_{s=1}^{b_n-1}\sum_{l=1}^{b_n-s}\frac{1}{b_n}\left(\sum_{i=1}^{l-1}U_iU_{s+i} + \sum_{i=n-b_n+l+1}^{n-s}U_iU_{s+i}\right)\right| \\
&\qquad + \frac{1}{n}\left|2\sum_{s=1}^{b_n-1}\sum_{l=1}^{b_n-s}\frac{\bar{B}_n}{b_n}\left(\sum_{i=1}^{l-1}(-U_i-U_{s+i})\right.\right. \\
&\qquad\qquad \left.\left. + \sum_{i=n-b_n+l+1}^{n-s}(-U_i-U_{s+i})\right)\right| \\
&\qquad + \frac{1}{n}\left|4\sum_{s=1}^{b_n-1}\sum_{l=1}^{b_n-s}\frac{\bar{B}_n^2}{b_n}(b_n-1)\right|.
\end{aligned}$$

We show that each of the three terms above tend to zero.

1. Letting

$$A(b_n) = \frac{2}{nb_n}\sum_{s=1}^{b_n-1}\sum_{l=1}^{b_n-s}\left(\sum_{i=1}^{l-1}U_iU_{s+i} + \sum_{i=n-b_n+l+1}^{n-s}U_iU_{s+i}\right),$$

it is straightforward to see that $A(b_n)$ can be written in the form $\sum_j\sum_k a_{jk} \times U_jU_k$ where the coefficients satisfy

$$0 \le a_{jk} \le \frac{2(b_n-1)}{nb_n}.$$

Since $U_1,\dots,U_{b_n-1},U_{n-b_n+2},\dots,U_n$ are i.i.d. N$(0,1)$ we have $EU_iU_{s+i}=0$, and hence $EA(b_n)=0$. Thus, we can apply Lemma 10 to obtain for $c\ge 1$ and a constant $K_c<\infty$,

$$EA(b_n)^{2c}\le K_c\left(\sum_j\sum_k a_{jk}^2\right)^c\le K_c\left(\frac{4(b_n-1)}{n}\right)^{2c}.$$

Then, similar to the argument given above, Markov's inequality and Borel–Cantelli yields $|A(b_n)|\to 0$ as $n\to\infty$ since by assumption there exists a constant $c\ge 1$ such that $\sum_n(b_n/n)^{2c}<\infty$.

2. Notice

$$\begin{aligned}
&\frac{2}{nb_n}\left|\sum_{s=1}^{b_n-1}\sum_{l=1}^{b_n-s}\bar{B}_n\left(\sum_{i=1}^{l-1}(-U_i-U_{s+i})+\sum_{i=n-b_n+l+1}^{n-s}(-U_i-U_{s+i})\right)\right|\\
&\quad\le\frac{2}{nb_n}|\bar{B}(n)|\sum_{s=1}^{b_n-1}\sum_{l=1}^{b_n-s}\left(\sum_{i=1}^{l-1}(|U_i|+|U_{s+i}|)\right.\\
&\qquad\qquad\left.+\sum_{i=n-b_n+l+1}^{n-s}(|U_i|+|U_{s+i}|)\right)\\
&\quad\le\frac{2}{nb_n}|\bar{B}(n)|\sum_{s=1}^{b_n-1}\sum_{l=1}^{b_n-s}\left(2\sum_{i=1}^{b_n}|U_i|+2\sum_{i=n-b_n+1}^{n}|U_i|\right)\\
&\quad=\frac{2(b_n-1)}{n^2}|B(n)|\left(\sum_{i=1}^{b_n}|U_i|+\sum_{i=n-b_n+1}^{n}|U_i|\right).
\end{aligned}$$

Then for $n$ large enough and $\epsilon>0$, Lemma 1 implies the right-hand side is bounded above by

$$2^{3/2}(1+\epsilon)\left[\frac{(b_n-1)^2\log\log n}{n^3}\right]^{1/2}\left(\sum_{i=1}^{b_n}|U_i|+\sum_{i=n-b_n+1}^{n}|U_i|\right).$$

Since $b_n^4n^{-3}\log\log n\to 0$ as $n\to\infty$ by assumption it suffices to show that

$$\frac{1}{b_n}\left(\sum_{i=1}^{b_n}|U_i|+\sum_{i=n-b_n+1}^{n}|U_i|\right)$$

stays bounded with probability 1 as $n\to\infty$.

The classical SLLN implies $b_n^{-1}\sum_{i=1}^{b_n}|U_i|$ stays bounded w.p. 1. All that remains is to show that $b_n^{-1}\sum_{i=n-b_n+1}^{n}|U_i|$ stays bounded almost surely as $n\to\infty$. Since the half-normal distribution enjoys a moment generating function we can appeal to the classical strong invariance principle (Komlós, Major and Tusnády [28, 29], Major [32]) to obtain that for sufficiently large $n$ with probability 1 there is a constant $C'$ such that

$$\left|\sum_{i=1}^{n}|U_i| - n\sqrt{2/\pi} - (1-2/\pi)B(n)\right| \le C'\log n. \tag{28}$$

Now observe that

$$\begin{aligned}
\frac{1}{b_n}\sum_{i=n-b_n+1}^{n}|U_i| &= \frac{1}{b_n}\left|\sum_{i=1}^{n}|U_i| - \sum_{i=1}^{n-b_n}|U_i|\right| \\
&= \frac{1}{b_n}\left|\left(\sum_{i=1}^{n}|U_i| - n\sqrt{2/\pi} - (1-2/\pi)B(n)\right)\right. \\
&\qquad - \left(\sum_{i=1}^{n-b_n}|U_i| - (n-b_n)\sqrt{2/\pi} - (1-2/\pi)B(n-b_n)\right) \\
&\qquad \left. + (1-2/\pi)(B(n)-B(n-b_n)) + b_n\sqrt{2/\pi}\right|.
\end{aligned}$$

Hence, by (28) and Lemma 2, for sufficiently large $n$,

$$\begin{aligned}
&\frac{1}{b_n}\sum_{i=n-b_n+1}^{n}|U_i| \\
&\quad \le \frac{1}{b_n}\left(2C'\log n + (1+\epsilon)\left[2b_n\left(\log\frac{n}{b_n} + \log\log n\right)\right]^{1/2} + b_n\sqrt{2/\pi}\right) \\
&\quad = \sqrt{2/\pi} + 2C'b_n^{-1}\log n + O((b_n^{-1}\log n)^{1/2}).
\end{aligned}$$

Hence, $b_n^{-1}\sum_{i=n-b_n+1}^{n}|U_i|$ stays bounded w.p. 1 since, by assumption, $b_n^{-1}\log n$ is bounded for all sufficiently large $n$.

3. Using Lemma 1 we have for every $\epsilon > 0$ and sufficiently large $n$

$$\begin{aligned}
\frac{4}{n}\left|\sum_{s=1}^{b_n-1}\sum_{l=1}^{b_n-s}\frac{\bar{B}_n^2}{b_n}(b_n-1)\right| &= \frac{2(b_n-1)^2}{n}\bar{B}_n^2 \\
&\le 4(1+\epsilon)^2\frac{(b_n-1)^2[\log\log n]}{n^2} \\
&\to 0
\end{aligned}$$

as $n\to\infty$ since we assumed $b_n^2 n^{-2}\log\log n \to 0$ as $n\to\infty$. $\square$

## APPENDIX C: MEAN-SQUARE CONSISTENCY PROOFS

**C.1. Preliminaries.** Recall that $\hat{\sigma}^2_{\mathrm{BM}}$ is defined at (6) while $\hat{\sigma}^2_{\mathrm{OBM}}$ is defined at (7). Next define the Brownian motion version of BM by

$$\tilde{\sigma}^2_{\mathrm{BM}} = \frac{b_n}{a_n - 1} \sum_{k=0}^{a_n - 1} (\bar{B}_k - \bar{B}_n)^2,$$

where $\bar{B}_k := b_n^{-1}(B((k+1)b_n) - B(kb_n + 1))$ for $k = 0, \ldots, a_n - 1$. Further define the Brownian motion version of OBM by

$$\tilde{\sigma}^2_{\mathrm{OBM}} = \frac{nb_n}{(n - b_n)(n - b_n + 1)} \sum_{j=0}^{n - b_n} (\bar{B}_j(b_n) - \bar{B}_n)^2.$$

LEMMA 11 (Damerdji [11]). *Suppose Assumption 2 holds, then*

$$E[\tilde{\sigma}^2_{\mathrm{OBM}}] = E[\tilde{\sigma}^2_{\mathrm{BM}}] = 1,$$

$$\frac{n}{b_n} \operatorname{Var}[\tilde{\sigma}^2_{\mathrm{OBM}}] = \frac{4}{3} + o(1)$$

*and*

$$\frac{n}{b_n} \operatorname{Var}[\tilde{\sigma}^2_{\mathrm{BM}}] = 2 + o(1),$$

*where the limits are taken as $n \to \infty$.*

The next claim follows from a careful examination of the proof of Lemma B.4 in Jones et al. [24].

LEMMA 12. *Suppose (13) holds with $\gamma(n) = n^\alpha \log n$ where $\alpha = 1/(2+\delta)$ and Assumption 2 holds. Then for sufficiently large $n$, there exist functions $g_1 : \mathbb{Z}_+ \to \mathbb{R}_+$ and $g_2 : \mathbb{Z}_+ \to \mathbb{R}_+$ such that*

$$|\hat{\sigma}^2_{\mathrm{BM}} - \sigma_g^2 \tilde{\sigma}^2_{\mathrm{BM}}| \le C^2 g_1(n) + C g_2(n),$$

*where the random variable $C$ is defined at (14). Moreover, if, as $n \to \infty$, $b_n^{-1} n^{2\alpha} [\log n]^3 \to 0$, then $g_1(n) \to 0$ and $g_2(n) \to 0$.*

We also require an analogous result for OBM.

LEMMA 13. *Suppose (13) holds with $\gamma(n) = n^\alpha \log n$ where $\alpha = 1/(2+\delta)$ and Assumption 2 holds. Then for sufficiently large $n$, there exist functions $g_1 : \mathbb{Z}_+ \to \mathbb{R}_+$, $g_2 : \mathbb{Z}_+ \to \mathbb{R}_+$, and $g_3 : \mathbb{Z}_+ \to \mathbb{R}_+$ such that*

$$|\hat{\sigma}^2_{\mathrm{OBM}} - \sigma_g^2 \tilde{\sigma}^2_{\mathrm{OBM}}| \le C^2 g_1(n) + C g_2(n) + (\hat{\sigma}^2_{\mathrm{OBM}} + \sigma_g^2 \tilde{\sigma}^2_{\mathrm{OBM}}) g_3(n),$$

*where the random variable $C$ is defined at (14). Moreover, if, as $n \to \infty$, $b_n^{-1} n^{2\alpha} [\log n]^3 \to 0$, then $g_1(n) \to 0$, $g_2(n) \to 0$, and $g_3(n) \to 0$.*

PROOF. Let $w$ be the modified Bartlett lag window. Notice

$$|\hat{\sigma}^2_{\mathrm{OBM}} - \sigma^2_g\tilde{\sigma}^2_{\mathrm{OBM}}| \le |\hat{\sigma}^2_{w,n} - \sigma^2_g\tilde{\sigma}^2_*| + |\hat{\sigma}^2_{\mathrm{OBM}} - \hat{\sigma}^2_{w,n}| + |\sigma^2_g\tilde{\sigma}^2_* - \sigma^2_g\tilde{\sigma}^2_{\mathrm{OBM}}|, \tag{29}$$

where $\hat{\sigma}^2_{w,n}$ and $\tilde{\sigma}^2_*$ are defined at (24) and (25), respectively. Our assumptions imply (15) and (16). Careful examination of the proof of Lemma 4 shows that for sufficiently large $n$,

$$|\hat{\sigma}^2_{w,n} - \sigma^2_g\tilde{\sigma}^2_*| \le C^2 g_1(n) + C g_2(n).$$

Let $g_3(n) = 1 - (n - b_n)(n - b_n + 1)/n^2$ and observe that

$$\begin{aligned}
|\hat{\sigma}^2_{\mathrm{OBM}} - \hat{\sigma}^2_{w,n}| &= \left| \left( \frac{n b_n}{(n - b_n)(n - b_n + 1)} - \frac{b_n}{n} \right) \sum_{j=0}^{n-b_n} [\bar{Y}_j(b_n) - \bar{Y}_n]^2 \right| \\
&= \frac{n b_n g_3(n)}{(n - b_n)(n - b_n + 1)} \sum_{j=0}^{n-b_n} [\bar{Y}_j(b_n) - \bar{Y}_n]^2 \\
&= g_3(n)\hat{\sigma}^2_{\mathrm{OBM}},
\end{aligned}$$

and, similarly, $|\sigma^2_g\tilde{\sigma}^2_* - \sigma^2_g\tilde{\sigma}^2_{\mathrm{OBM}}| = \sigma^2_g\tilde{\sigma}^2_{\mathrm{OBM}} g_3(n)$. Combining these results with (29) yields the claim. □

LEMMA 14. *Let $\hat{\sigma}^2_n$ be either the BM or OBM estimator of $\sigma^2_g$ and suppose (13) holds with $\gamma(n) = n^\alpha \log n$ where $\alpha = 1/(2+\delta)$ and Assumption 2 holds. Also, assume that $EC^4 < \infty$ where $C$ is defined at (14) and $E_\pi g^4 < \infty$. If $b_n^{-1} n^{2\alpha}[\log n]^3 \to 0$ as $n \to \infty$, then $E[|\hat{\sigma}^2_n - \sigma^2_g\tilde{\sigma}^2|] \to 0$ and $E[(\hat{\sigma}^2_n - \sigma^2_g\tilde{\sigma}^2)^2] \to 0$ while if as $n \to \infty$, $b_n^{-1} n^{1/2+\alpha}[\log n]^{3/2} \to 0$, then $E[\frac{n}{b_n}(\hat{\sigma}^2_n - \sigma^2_g\tilde{\sigma}^2)^2] \to 0$.*

PROOF. We will prove only the first claim for BM as the other proofs are quite similar. The omitted proofs for OBM require the use of Lemma 13 in place of Lemma 12 and Lemma 4 with the ergodic theorem in place of Lemma B.4 from Jones et al. [24].

From Lemma 12 there exists an integer $N_0$ and functions $g_1$ and $g_2$ such that

$$\begin{aligned}
|\hat{\sigma}^2_{\mathrm{BM}} - \sigma^2_g\tilde{\sigma}^2_{\mathrm{BM}}| &= |\hat{\sigma}^2_{\mathrm{BM}} - \sigma^2_g\tilde{\sigma}^2_{\mathrm{BM}}| I(0 \le n \le N_0) + |\hat{\sigma}^2_{\mathrm{BM}} - \sigma^2_g\tilde{\sigma}^2_{\mathrm{BM}}| I(N_0 < n) \\
&\le |\hat{\sigma}^2_{\mathrm{BM}} - \sigma^2_g\tilde{\sigma}^2_{\mathrm{BM}}| I(0 \le n \le N_0) \\
&\quad + [C^2 g_1(n) + C g_2(n)] I(N_0 < n) \\
&:= g_n(X_0, \dots, X_n, B(0), \dots, B(n)).
\end{aligned}$$

Now

$$\begin{aligned}&E|g_n(X_0,\dots,X_n,B(0),\dots,B(n))|\\&\qquad\le E|\hat\sigma^2_{\mathrm{BM}}-\sigma_g^2\tilde\sigma^2_{\mathrm{BM}}|+[g_1(n)E(C^2)+g_2(n)E(C)],\end{aligned}$$

and, since Lemma 11 and our assumptions on the moments of $g$ imply

$$E|\hat\sigma^2_{\mathrm{BM}}-\sigma_g^2\tilde\sigma^2_{\mathrm{BM}}|\le E\hat\sigma^2_{\mathrm{BM}}+\sigma_g^2E\tilde\sigma^2_{\mathrm{BM}}=E\hat\sigma^2_{\mathrm{BM}}+\sigma_g^2<\infty,$$

it follows from our assumptions on the moments of $C$ that $E|g_n|<\infty$. Next observe that as $n\to\infty$, we have $g_n\to0$ w.p. 1 and $Eg_n\to0$ by Lemma 12. From Lemma B.4 in Jones et al. [24] we have that $|\hat\sigma^2_{\mathrm{BM}}-\sigma_g^2\tilde\sigma^2_{\mathrm{BM}}|\to0$ w.p. 1 as $n\to\infty$. An application of the generalized majorized convergence theorem (Zeidler [54], page 1015) implies that, as $n\to\infty$, $E[|\hat\sigma^2_{\mathrm{BM}}-\sigma_g^2\tilde\sigma^2_{\mathrm{BM}}|]\to0$. □

**C.2. Proof of Theorem 3.** We will prove only the claim for BM as the proof for OBM is nearly identical.

Recall the mean-square error (MSE) of the estimator $\hat\sigma^2_{\mathrm{BM}}$ of $\sigma_g^2$ is given by

$$\mathrm{MSE}(\hat\sigma^2_{\mathrm{BM}})=\mathrm{Var}(\hat\sigma^2_{\mathrm{BM}})+\mathrm{Bias}^2(\hat\sigma^2_{\mathrm{BM}}).$$

First, consider the bias term. Recall from Lemma 11 that $E\tilde\sigma^2_{\mathrm{BM}}=1$ so that

$$\mathrm{Bias}(\hat\sigma^2_{\mathrm{BM}})=E(\hat\sigma^2_{\mathrm{BM}})-\sigma_g^2\le E[|\hat\sigma^2_{\mathrm{BM}}-\sigma_g^2\tilde\sigma^2_{\mathrm{BM}}|]\to0\qquad\text{as } n\to\infty,$$

by Lemma 14. Hence $\mathrm{Bias}^2(\hat\sigma^2_{\mathrm{BM}})=o(1)$ as $n\to\infty$.

Note that the claim will follow if we can show that $\mathrm{Var}(\hat\sigma^2_{\mathrm{BM}})=o(1)$ as $n\to\infty$. Thus we now consider the variance term but we begin with a preliminary result. Recall from Lemma 11 that

$$\frac{n}{b_n}\mathrm{Var}[\tilde\sigma^2_{\mathrm{BM}}]=2+o(1).$$

Define

$$\begin{aligned}\eta(\hat\sigma^2_{\mathrm{BM}},\tilde\sigma^2_{\mathrm{BM}})&=\mathrm{Var}(\hat\sigma^2_{\mathrm{BM}}-\sigma_g^2\tilde\sigma^2_{\mathrm{BM}})\\&\qquad+2\sigma_g^2E[(\tilde\sigma^2_{\mathrm{BM}}-E\tilde\sigma^2_{\mathrm{BM}})(\hat\sigma^2_{\mathrm{BM}}-\sigma_g^2\tilde\sigma^2_{\mathrm{BM}})].\end{aligned}\tag{30}$$

Using the Cauchy–Schwarz inequality and the the fact that $\mathrm{Var}(X)\le EX^2$ obtain

$$\begin{aligned}|\eta|&=|\mathrm{Var}(\hat\sigma^2_{\mathrm{BM}}-\sigma_g^2\tilde\sigma^2_{\mathrm{BM}})+2\sigma_g^2E[(\tilde\sigma^2_{\mathrm{BM}}-E\tilde\sigma^2_{\mathrm{BM}})(\hat\sigma^2_{\mathrm{BM}}-\sigma_g^2\tilde\sigma^2_{\mathrm{BM}})]|\\&\le E[(\hat\sigma^2_{\mathrm{BM}}-\sigma_g^2\tilde\sigma^2_{\mathrm{BM}})^2]+2\sigma_g^2(E[(\hat\sigma^2_{\mathrm{BM}}-\sigma_g^2\tilde\sigma^2_{\mathrm{BM}})^2]\mathrm{Var}(\tilde\sigma^2_{\mathrm{BM}}))^{1/2}\\&=o(1)+2\sigma_g^2\left(\frac{b_n}{n}\right)^{1/2}[o(1)(2+o(1))]^{1/2}\qquad\text{as } n\to\infty,\end{aligned}$$

and hence $\eta(\hat{\sigma}^2_{\mathrm{BM}}, \tilde{\sigma}^2_{\mathrm{BM}}) = o(1)$ since $b_n/n \to 0$ as $n \to \infty$. Now

$$\begin{aligned}
\mathrm{Var}(\hat{\sigma}^2_{\mathrm{BM}}) &= E[((\hat{\sigma}^2_{\mathrm{BM}} - \sigma_g^2\tilde{\sigma}^2_{\mathrm{BM}}) \\
&\qquad + \sigma_g^2(\tilde{\sigma}^2_{\mathrm{BM}} - E\tilde{\sigma}^2_{\mathrm{BM}}) - (E\hat{\sigma}^2_{\mathrm{BM}} - \sigma_g^2 E\tilde{\sigma}^2_{\mathrm{BM}}))^2] \\
&= \sigma_g^4 E[(\tilde{\sigma}^2_{\mathrm{BM}} - E\tilde{\sigma}^2_{\mathrm{BM}})^2] \\
&\quad + E[((\hat{\sigma}^2_{\mathrm{BM}} - \sigma_g^2\tilde{\sigma}^2_{\mathrm{BM}}) - E(\hat{\sigma}^2_{\mathrm{BM}} - \sigma_g^2\tilde{\sigma}^2_{\mathrm{BM}}))^2] \\
&\quad + 2\sigma_g^2 E[(\tilde{\sigma}^2_{\mathrm{BM}} - E\tilde{\sigma}^2_{\mathrm{BM}})((\hat{\sigma}^2_{\mathrm{BM}} - \sigma_g^2\tilde{\sigma}^2_{\mathrm{BM}}) \\
&\qquad\qquad - E(\hat{\sigma}^2_{\mathrm{BM}} - \sigma_g^2\tilde{\sigma}^2_{\mathrm{BM}}))] \\
&= \sigma_g^4 \mathrm{Var}(\tilde{\sigma}^2_{\mathrm{BM}}) + \eta(\hat{\sigma}^2_{\mathrm{BM}}, \tilde{\sigma}^2_{\mathrm{BM}}) \\
&= 2\sigma_g^4 \frac{b_n}{n} + o\left(\frac{b_n}{n}\right) + o(1) \qquad \text{as } n \to \infty.
\end{aligned} \tag{31}$$

Therefore, we conclude that $\mathrm{Var}(\hat{\sigma}^2_{\mathrm{BM}}) = o(1)$ as $n \to \infty$.

**C.3. Proof of Theorem 4.** We only prove the claim for BM as the proof for OBM is nearly identical.

From (31) and Lemma 11 we obtain

$$\frac{n}{b_n}\mathrm{Var}(\hat{\sigma}^2_{\mathrm{BM}}) = 2\sigma_g^4 + o(1) + \frac{n}{b_n}\eta,$$

where $\eta$ is defined at (30). Note that the claim will follow if we show that $\frac{n}{b_n}\eta = o(1)$ as $n \to \infty$. As in the proof of Theorem 3 we have

$$\begin{aligned}
\frac{n}{b_n}|\eta| &= \frac{n}{b_n}|\mathrm{Var}(\hat{\sigma}^2_{\mathrm{BM}} - \sigma_g^2\tilde{\sigma}^2_{\mathrm{BM}}) + 2\sigma_g^2 E[(\tilde{\sigma}^2_{\mathrm{BM}} - E\tilde{\sigma}^2_{\mathrm{BM}})(\hat{\sigma}^2_{\mathrm{BM}} - \sigma_g^2\tilde{\sigma}^2_{\mathrm{BM}})]| \\
&\le E\left[\frac{n}{b_n}(\hat{\sigma}^2_{\mathrm{BM}} - \sigma_g^2\tilde{\sigma}^2_{\mathrm{BM}})^2\right] \\
&\quad + 2\sigma_g^2\left(E\left[\frac{n}{b_n}(\hat{\sigma}^2_{\mathrm{BM}} - \sigma_g^2\tilde{\sigma}^2_{\mathrm{BM}})^2\right]\frac{n}{b_n}\mathrm{Var}(\tilde{\sigma}^2_{\mathrm{BM}})\right)^{1/2} \\
&= o(1) + 2\sigma_g^2(o(1)(2 + o(1)))^{1/2} \qquad \text{as } n \to \infty,
\end{aligned}$$

using the results of Lemmas 11 and 14. Hence $\frac{n}{b_n}\eta = o(1)$ as $n \to \infty$.

**Acknowledgments.** The authors are grateful to two referees and the associate editor for their constructive comments which resulted in many improvements.

Department of Statistics
University of California
Riverside, California 92521
USA
E-mail: jflegal@ucr.edu

School of Statistics
University of Minnesota
Minneapolis, Minnesota 55455
USA
E-mail: galin@stat.umn.edu